\documentclass[12pt,a4paper,twoside,leqno]{article}
\usepackage{amsthm,amsmath,amssymb,mathabx}
\usepackage{color}
\usepackage{graphicx}
\usepackage{graphicx}
\usepackage{setspace}
\usepackage{esint}
\usepackage{epstopdf}

\usepackage{placeins}

\usepackage{setspace}

\newtheorem{teor}{Theorem}[section]
\newtheorem{defin}[teor]{Definition}
\newtheorem{lemm}[teor]{Lemma}
\newtheorem{osse}[teor]{Remark}
\newtheorem{prop}[teor]{Proposition}
\newtheorem{defi}[teor]{Definition}
\newtheorem{coro}[teor]{Corollary}
\newtheorem{prob}[teor]{Problem}
\newtheorem{hypo}[teor]{Hypothesis}

\newcommand{\bele}{\begin{lemm}\begin{sl}}
\newcommand{\enle}{\end{sl}\end{lemm}}
\newcommand{\bedef}{\begin{defi}\begin{sl}}
\newcommand{\eddef}{\end{sl}\end{defi}}
\newcommand{\bete}{\begin{teor}\begin{sl}}
\newcommand{\ente}{\end{sl}\end{teor}}
\newcommand{\beos}{\begin{osse}\begin{sl}}
\newcommand{\eddos}{\end{sl}\end{osse}}
\newcommand{\bepr}{\begin{prop}\begin{sl}}
\newcommand{\empr}{\end{sl}\end{prop}}
\newcommand{\bepro}{\begin{prob}\begin{rm}}
\newcommand{\empro}{\end{rm}\end{prob}}
\newcommand{\bede}{\begin{defin}\begin{sl}}
\newcommand{\edde}{\end{sl}\end{defin}}
\newcommand{\beco}{\begin{coro}\begin{sl}}
\newcommand{\enco}{\end{sl}\end{coro}}
\newcommand{\behy}{\begin{hypo}\begin{sl}}
\newcommand{\enhy}{\end{sl}\end{hypo}}

\newcommand{\thspace}{\hspace{3mm}}

\newcommand{\beeq}[1]{\begin{equation}\label{#1}}
\newcommand{\eddeq}{\end{equation}}

\newcommand{\beeqa}[1]{\begin{eqnarray}\label{#1}}
\newcommand{\eddeqa}{\end{eqnarray}}

\newcommand{\beal}[1]{\begin{align}\label{#1}}
\newcommand{\eddal}{\end{align}}

\newcommand{\bespl}[1]{\begin{split}\label{#1}}
\newcommand{\edspl}{\end{split}}

\newcommand{\bega}[1]{\begin{gather}\label{#1}}
\newcommand{\edga}{\end{gather}}

\newcommand{\beeqax}{\begin{eqnarray*}}
\newcommand{\eddeqax}{\end{eqnarray*}}

\def\qed{\ifmmode 
  \else \leavevmode\unskip\penalty9999 \hbox{}\nobreak\hfill
  \fi
  \quad\hbox{\hskip.5em\vrule width.4em height.6em depth.05em\hskip.1em}}
\def\endproofsym{\qed}

\def\endnobox{\def\endproofsym{}\end{proof}\def\endproofsym{\qed}}

\newcommand{\no}{\nonumber}

\newcommand{\beeqao}{\begin{eqnarray}\no}
\newcommand{\bealo}{\begin{align}\no}
\newcommand{\besplo}{\begin{split}\no}
\newcommand{\begao}{\begin{gather}\no}

\newcommand{\elefin}{\color{black}}

\newcommand{\EEE}{\color{black}}

\newcommand{\ele}{\color{black}}

\newcommand{\elee}{\color{black}}

\numberwithin{equation}{section}
\begin{document}

\title{A nonlinear model for marble sulphation including surface rugosity: theoretical and numerical results}

\author{
  Elena Bonetti
  \footnote{Dipartimento di Matematica ``F. Enriques'',
Universit\`a degli Studi di Milano, Via C. Saldini 50, 20133 Milano,
Italy. Istituto di Matematica Applicata e Tecnologie Informatiche ``Enrico Magenes'', CNR, Via Ferrata 1, 27100 Pavia, Italy. \texttt{elena.bonetti@unimi.it}
 }
  \and
  Cecilia Cavaterra
  \footnote{Dipartimento di Matematica ``F. Enriques'',
Universit\`a degli Studi di Milano, Via C. Saldini 50, 20133 Milano,
Italy. Istituto di Matematica Applicata e Tecnologie Informatiche ``Enrico Magenes'', CNR, Via Ferrata 1, 27100 Pavia, Italy. \texttt{cecilia.cavaterra@unimi.it}
 }
 \and
  Francesco Freddi
  \footnote{Dipartimento di Ingegneria e Architettura, Universit\`a degli Studi di Parma, Parco Area delle Scienze, 181/A, 43124 Parma, Italy.
  \texttt{francesco.freddi@unipr.it}
  }
  \and
  Maurizio Grasselli
  \footnote{Dipartimento di Matematica, Politecnico di Milano, Via E. Bonardi 9, 20133 Milano, Italy.
    \texttt{maurizio.grasselli@polimi.it}}
    \and
  Roberto Natalini
  \footnote{Istituto per le Applicazioni del Calcolo ``M. Picone'', CNR, Via dei Taurini 19, 00185 Roma, Italy.
    \texttt{roberto.natalini@cnr.it}}}

\maketitle

\begin{abstract}
\noindent
We consider an evolution system describing the phenomenon of marble sulphation of a monument, accounting of the surface rugosity.
We first prove a local in time well posedness result. Then, stronger assumptions on the data allow us to establish the existence of a global in
time solution. Finally, we perform some numerical simulations that illustrate the main feature of the proposed model.
\end{abstract}

\vspace{.2cm}


\vspace{2mm}

\noindent
{\bf AMS (MOS) subject clas\-si\-fi\-ca\-tion:}\thspace 35M33, 65M06, 76R50.

\section{Introduction}

Deterioration and damage of stones is a complex problem and one of the main concern for people
working in the field of conservation and restoration of cultural heritage. It is extremely difficult
to isolate a single factor in this kind of processes, which are the results of the interaction of various
mechanisms. These processes can be described through free boundary models, as in the case of damage induced
by pollution, or by using a phase field approach.
We introduce a new differential system coupling bulk and surface evolution equations to describe the
phenomenon of marble sulphation of a monument, including the surface rugosity.

More precisely,  we consider a monument made of calcium carbonate stone, located in a smooth bounded domain
$\Omega\subseteq\mathbb{R}^3$, with boundary $\Gamma$, and subjected to a degradation process during a time
interval $(0,T)$, for any given $T>0$.

Following the approach introduced in \cite{N1, N2, N3, GN1, GN3D}  for related models of marble sulphation,  we consider
the system
\begin{align}\label{reazdiff}
 &\partial_t(\phi(c)s)-\hbox{div }(\phi(c)\nabla s)=-\lambda\phi(c)sc,\quad\hbox{in }\Omega \times (0,T),\\\label{reazchim}
 &\partial_t c=-\lambda\phi(c)sc,\quad\hbox{in }\Omega \times (0,T),\\\label{boundflux}
&\phi(c)\partial_ns=-\nu(r)(s -\bar s),\quad\hbox{on }\Gamma\times (0,T),\\\label{surfacedam}
&\partial_t r+\partial W(r)+\Psi'(r)+G(r,c,s)\ni F, \quad\hbox{on }\Gamma \times (0,T),
\end{align}
equipped with the set of initial  conditions
\begin{equation}\label{cauchy}
s(0)=s_0,\, \hbox{ in }\, \Omega, \quad c(0)=c_0 \, \hbox{ in }\, \Omega , \quad r(0)=r_0 \, \hbox{ on }\, \Gamma.
\end{equation}

Here $s$ stands for the $SO_2$ porous concentration inside the material,  and $c$ for the local density of $CaCO_3$.
On the boundary $\Gamma$ the variable $r\geq 0$ denotes the rugosity of the surface.

The first equation describes the evolution of the porous concentration of $SO_2$, where $\phi(c)$ is the porosity of the medium.
This evolution is driven by the Fick's law and on the right hand-side we have the reaction term between  $SO_2$  and calcium carbonate,
with rate $\lambda >0$.
The second equation takes into account for the loss of calcium carbonate due the reaction.

The novelty in this model with respect to \cite{N1,N3} is the introduction of the Robin boundary condition \eqref{boundflux} for the flux of $SO_2$
through the external boundary,depending not only on the porosity of the medium and the external concentration $\bar s$, but also on an effective
permeability coefficient $\nu$, which is a function of the superficial rugosity $r$.

Let us now describe what $r$ represents. From the physical point of view the rugosity is quite a complex quantity, corresponding locally
to the microscopic variation of a surface with respect to a flat configuration. For a precise definition see, for instance, \cite{whitehouse,iso}.
Actually, in our model $r$ stands for an effective macroscopic parameter which is formally a local damage parameter. As $r$ grows, the profile of the
surface has a greater deviation from the flat configuration and the microscopic density of peaks and valleys increases.
In practice, this corresponds to a greater surface exposed to the pollution action, implying an increasing of the $SO_2$ flux permeability coefficient.
Here,  the value $r=0$ corresponds to a completely smooth surface. Moreover, we deal with a material in which the ``direction'' of the  evolution of $r$
is not a priori fixed, i.e. $r$ may increase or decrease through different external actions or repairing itself.
The evolution of $r$ is governed by the damage-like equation \eqref{surfacedam} (cf., e.g., \cite{Bonetti,Fremond,Natalini-BARENBLAT} and
the references therein).
The symbol $\partial W$ stands for the subdifferential of a proper, convex, and lower semicontinuous function  $W$, accounting for possible internal
constraints on $r$.  The function $G$ depends in particular on $c$ and $s$ through the porosity function $\phi(c)$ and satisfies  $G(r,0,s)=G(r,c,0)=0$.
The function $\Psi'$ is sufficiently smooth and accounts for the non-monotone dynamics of $r$. The function $F$ on the right hand side of
\eqref{surfacedam} denotes possible external actions responsible for the formation of rugosity on the boundary, like, e.g., wind, rain or temperature variations.

It is worth noting that, as many dissipative evolution equations, our system  \eqref{reazdiff}-\eqref{surfacedam} has a gradient flow structure
(see, e.g., \cite{peletier}). Indeed, let us introduce the energy functionals
\begin{align}
&{\cal E}_1[c,s,r]=\int_\Omega\left (\frac {\phi(c)}2|\nabla s|^2+\lambda c\phi(c)\frac {s^2}2-\lambda\phi'(c)\phi(c)c\frac {s^3}3\right)
+\int_\Gamma\frac {\nu(r)}2|s-\bar s|^2,\\\no
&{\cal E}_2[c,s,r]=\int_\Gamma (W(r)+\Psi(r)+\widehat G(r,c,s)-Fr),
\end{align}
where  $\widehat G$ is such that $G=\partial_r\widehat G$.
Then, using the kinetic equation \eqref{reazchim} and letting
\begin{align}
&{\cal R}_1[\partial_t s]=\frac 1 2\int_\Omega\phi(c)|\partial_t s|^2,\\
&{\cal R}_2[\partial_t r]=\frac 1 2\int_\Gamma |\partial_t r|^2,
\end{align}
we can formally write \eqref{reazdiff}, combined with \eqref{boundflux} and \eqref{surfacedam}, as
 \begin{align}\label{GFreaz}
&\partial{\cal R}_1+\partial_s{\cal E}_1=0,\\
&\partial{\cal R}_2+\partial_r{\cal E}_2\ni0.
 \end{align}
Actually here we do not exploit this structure since we use only standard energy arguments, but this remark could be helpful for future developments.

In this paper we first prove a local in time well posedness result for system \eqref{reazdiff}--\eqref{cauchy} in finite energy spaces.
Next, under slightly stronger assumptions, we can establish some uniform bounds for the solution which in turn imply the global in time existence.
Finally, we perform some numerical finite element simulations to assess the behavior of our model according to different physical situations.

The plan of the paper is the following. Section 2 is devoted to the statement of the main theorems and assumptions. The proof of the local existence
result is detailed in Section 3. Section 4 contains the proof of the global existence result. Finally, in Section 5 we present the numerical simulations.


 \section{Main  results}

In this section, we introduce the formulation of system \eqref{reazdiff}--\eqref{cauchy} in the correct functional setting and we state the main results.
From now on, for the sake of simplicity, we will take $\lambda=1$.

First of all,  we assume that the porosity function  $\phi$ is a linear function of the density $c$, i.e.,
\begin{equation}\label{sceltafi}
\phi(c)=A+Bc,
\end{equation}
where $A$ and $B$ are constants (see \cite{GN1} and the references therein).

Concerning the data, throughout the paper we make use of the following assumptions
\begin{align}
\label{A7}&\ele W:[0,+\infty)\rightarrow[0,+\infty]\hbox{ is convex and l.s.c.},\,W(0)=0,\\\no
&\Psi\in W^{2,\infty}({\mathbb R}), \quad G\in W^{1,\infty}({\mathbb R}^3), \quad F\in L^2(0,T;L^2(\Gamma)),\\
\label{A3}&\bar s\in H^{1/2}(\Gamma),\quad \bar s\geq 0, \, \hbox{ a.e. in }\Gamma,\\
\label{A2}&\nu\in W^{1,\infty}(\mathbb{R}),\quad\nu \geq 0, \, \hbox{ a.e. in }\mathbb{R}, \\
\label{A4}&c_0\in H^2(\Omega), \quad 0\leq c_0(x)\leq C_0, \, \, \, \forall \, x \in \Omega,\\
\label{A6}&s_0\in H^2(\Omega), \quad s_0(x) \geq 0, \,\,\, \forall \, x \in \Omega, \\
\label{A5}&r_0\in L^2(\Gamma), \quad W(r_0)\in L^1(\Gamma), \quad r_0 \geq 0, \, \hbox{ a.e. in } \Gamma,\\
\label{A5bis}&\phi(c_0)\partial_ns_0 = - \nu(r_0)(s_0 - \bar s), \, \hbox{ a.e. in }\Gamma,\\
\label{A1}&A>0, \, \, A + B C_0 >0, \\
\label{A9}&B\leq \frac1S_0\,\hbox{ and }\,\bar s\leq S_0,\hbox{ a.e. in }\Gamma,\,\hbox{ with }S_0>0\,
\hbox{ s.t. }s_0(x)\leq S_0, \,\,\, \forall \, x \in \Omega.
\end{align}

\begin{osse}\label{stazs1}
We are considering the case in which $\bar s$ does not depend on $t$, i.e. $\partial_t\bar s=0$, for the sake
of simplicity. Note in addition that, by Sobolev's embeddings, $s\in L^4(\Gamma)$.
\end{osse}

\begin{osse}\label{stazs}
We recall that if $W$ is convex and lower semicontinuous then the subdifferential $\partial W$ is a maximal monotone operator defined as
$$\partial W (\sigma_0) = \{ \xi  \in \mathbb{R} \, : \, W(\sigma) \geq W(\sigma_0) + \xi(\sigma - \sigma_0)\}.$$
\end{osse}

\begin{osse}\label{stazs3}
It is reasonable to impose a constraint on $r$ so that $r\in[0,R_0]$, where $R_0>0$ depends on the crystals dimension for the material we are considering.
Actually, we can ensure the validity of this internal constraint assuming in \eqref{A7}, e.g., $W(x)=I_{[0,R_0]}(x)$.
\end{osse}


We can now formulate our problem

\vskip0.3truecm
\noindent{\bf Problem $(P)$}.\quad
Find a triplet $(s,c,r)$  such that
\begin{align}
&s\in H^1(0,T;L^2(\Omega))\cap L^\infty(0,T;H^1(\Omega))\cap L^2(0,T;H^2(\Omega)),\\
&c\in  H^1(0,T;H^2(\Omega)),\\
&r\in H^1(0,T;L^2(\Gamma)),\\
&s\geq 0, \quad  \, 0 \leq c \leq C_0, \quad \hbox{ a.e. in }\Omega\times(0,T),
\end{align}
and satisfying
\begin{align}\label{eq1w}
&\partial_t(\phi(c)s)-\hbox{div }(\phi(c)\nabla s)=-\phi(c)sc,  \quad\hbox{a.e. \!in } \Omega\times(0,T),\\\label{eq1wbis}
&\phi(c)\partial_ns=-\nu(r)(s-\bar s), \quad\hbox{a.e. \!on } \Gamma\times(0,T),\\\label{eq2w}
&\partial_t c=-\phi(c)cs,\quad\hbox{a.e. in }\Omega \times (0,T),\\\label{eq3w}
&\partial_t r+\xi+\Psi'(r)+G(r,c,s)=F,\quad \xi\in\partial W(r),\quad\hbox{a.e. on }\Gamma \times (0,T),\\\label{eq4w}
&s(0)=s_0, \,\,\, c(0)=c_0, \quad\hbox{ a.e. in }\, \Omega , \quad\quad r(0)=r_0, \quad \hbox{ a.e. on }\, \Gamma.
\end{align}

We can now state our main existence results. The first theorem concerns the existence of a local in time solution.

\begin{teor}\label{eloc}
Let $T>0$ and \eqref{A7}--\eqref{A1} hold. Then there exists a time $\widehat T \in (0,T]$ such that Problem $(P)$ admits a unique solution in $(0,\widehat T)$.
Moreover, the following properties hold
\begin{align}\label{loc1}
& s\in W^{1,\infty}(0,\widehat T;L^2(\Omega))\cap H^1(0,\widehat T;H^1(\Omega))\cap L^\infty(0,\widehat T;H^2(\Omega)),\\\label{loc2}
& c\in W^{1,\infty}(0,\widehat T;H^2(\Omega))\cap W^{2,\infty}(0,\widehat T;L^2(\Omega)),\\\label{loc3}
& s\geq 0, \quad  \, 0 \leq c \leq C_0, \quad \hbox{ a.e. in } \Omega\times (0,\widehat T),\\\label{loc4}
&\xi\in L^2(0,\widehat T;L^2(\Gamma)).
\end{align}
\end{teor}

Assumption  \eqref{A9} allows us to extend the existence of the solution to the whole time interval $(0,T)$.

\begin{teor}\label{eglob}
Let $T>0$ and \eqref{A7}--\eqref{A9} hold. Then there exists a unique global solution to Problem $(P)$ on the whole time interval $(0,T)$. Moreover, the
following properties hold
\begin{align}\label{glob1}
& s\in W^{1,\infty}(0,T;L^2(\Omega))\cap H^1(0,T;H^1(\Omega))\cap L^\infty(0,T;H^2(\Omega)), \\\label{glob2}
& c\in W^{1,\infty}(0,T;H^2(\Omega)),  \\\label{glob3}
& 0\leq s\leq S_0, \quad 0 \leq c \leq C_0, \quad \hbox{ a.e. in }\Omega\times (0,T),\\\label{glob4}
&\xi\in L^2(0,T;L^2(\Gamma)).
\end{align}
\end{teor}

\begin{osse}
Let us point out that, due to the uniform bound on $c$ and $s$ (see  \eqref{glob3}), our results also hold if  $G$  locally Lipschitz with respect
to $c$ and $s$ (see \eqref{A7}).
\end{osse}

\section{Proof of Theorem \ref{eloc}}

We consider problem \eqref{eq1w}--\eqref{eq4w}. Let us define, for $R>\|s_0\|_{L^2(\Omega)}$  and $T>0$, the subset

\begin{align}
{\cal X}_{R,T}:=\big\{&s\in C([0,T];L^2(\Omega)) \cap L^2(0,T;H^2(\Omega)) \, : \, \\\nonumber
&\|s\|_{ L^\infty(0,T;L^2(\Omega))\cap L^2(0,T;H^2(\Omega))}\leq R, \quad
s\geq 0, \quad \hbox{ a.e. in }\Omega\times(0,T)\big\},
\end{align}
where $T$ will be fixed later on as a function of $R$.
We aim to construct a mapping $${\cal S}:{\cal X}_{R,T}\rightarrow{\cal X}_{R,T}$$ which results to be a contraction for a suitable choice of $T$,
with respect to the norm of $C([0,T];L^2(\Omega))\cap L^2(0,T;H^2(\Omega))$.
As a consequence,  we will deduce that it admits a unique fixed point and we will show that this fixed point provides a solution to Problem $(P)$.

By abusing of notation, in the following we will use the same symbol $C$ for possibly different positive constants, depending only on the data of the problem,
on $\Omega$ and (continuously) on $T$ at most; while $C(R)$ will denote a positive constant which also depends on $R$.

\subsection{First step: finding $c$}

Let us fix $\widehat s \in{\cal X}_{R,T}$ and consider the following problem (cf. \eqref{reazchim}, \eqref{cauchy} and \eqref{sceltafi})
\begin{equation}\label{eqPF1}
\partial_t c=-Ac\widehat s-Bc^2\widehat s , \quad  c(0)=c_0 \quad\hbox{a.e. in }\Omega.
\end{equation}
This is actually a Cauchy problem (in time) for a Bernoulli type equation for $c$, which has the explicit solution
\begin{equation}\label{esplicitoc}
\widehat c(x,t)=\frac{Ac_0(x)}{(A+Bc_0(x))e^{A\int_0^t\widehat s(x,\tau)d\tau}-Bc_0(x)},\quad \hbox{ in } \Omega\times(0,T).
\end{equation}
Consequently, by \eqref{esplicitoc}, we can define an operator ${\cal S}_1(\widehat s)=\widehat c$.  Then, exploiting the bounds on $c_0$ (see \eqref{A4})
and on $A$ and $B$ (see \eqref{A1}), we deduce
\begin{equation}\label{stimaPF1ter}
A  + B c_0(x) \geq
\begin{cases}
&A > 0, \hskip1.4truecm \hbox{ if } B \geq0, \\
&A + BC_0 >0, \,\, \hbox{ if } B < 0,
\end{cases}
\quad \, \forall \, x \in \Omega.
\end{equation}
So that, we have
$$(A+Bc_0(x))e^{A\int_0^t\widehat s(x,\tau)d\tau}-Bc_0(x) > A+Bc_0(x) -Bc_0(x) >A>0, \, \hbox{ in } \Omega\times(0,T),
$$
and then (cf. also \eqref{esplicitoc})
\begin{equation}\label{stimaPF1}
0\leq \widehat c(x,t)\leq c_0(x) \leq C_0, \, \hbox{ in } \Omega\times(0,T),
\end{equation}
\begin{equation}
\partial_t\widehat c(x,t) \leq 0, \, \hbox{ a.e. in } \Omega\times(0,T).
\end{equation}
As a consequence we deduce, for all $(x,t) \in \Omega \times (0,T)$, the following inequalities
\begin{align}
&0 < A \leq A + B \widehat c(x,t)  \leq A + BC_0, \, \, \, \hbox{ if } B \geq 0, \\\no
&0 < A + BC_0 \leq A + B \widehat c(x,t) \leq A, \, \, \, \hbox{ if } B \geq 0.
\end{align}
Then there exist two positive constants $m$ and $M$ such that
\begin{equation}\label{stimaPF0}
0 < m \leq A + B\widehat c(x,t) = \phi(\widehat c) \leq M, \, \hbox{ in } \Omega\times(0,T).
\end{equation}
 Moreover, \eqref{esplicitoc} and the definition of {\ele ${\cal X}_{R,T}$} imply
\begin{equation}\label{stimaPF1bis}
\|\widehat c\|_{L^2(0,T;H^2(\Omega))}\leq C(R).
\end{equation}
By \eqref{eqPF1}
we also get  (see \eqref{sceltafi})
\begin{align}\label{stimaPF1tris}
&\|\partial_t \widehat c\|^2_{L^2(0,T;H^2(\Omega))}=\int_0^t\|\phi(\widehat c)\widehat c\widehat s\|^2_{H^2(\Omega)}\\\no
&\leq \int_0^t\left(\|\phi(\widehat c)\widehat c\|^2_{L^\infty(\Omega)}\|\widehat s\|^2_{H^2(\Omega)}
+ \|\phi(\widehat c)\widehat c\|^2_{H^2(\Omega)}\|\widehat s\|^2_{L^\infty(\Omega)}\right)\leq C(R).
\end{align}
Thus, we eventually obtain
\begin{equation}\label{stimaPF2}
\|\widehat c\|_{H^1(0,T;H^2(\Omega))}\leq C(R).
\end{equation}

\subsection{Second step: finding $r$}

We indicate by $\widehat s_{|_\Gamma}$  the trace on $\Gamma$ of the fixed $\widehat s\in {\cal X}_R$. Analogously, $\widehat c_{|_{\Gamma}}$ indicates the trace
of  $\widehat c$ (given by \eqref{esplicitoc}) on $\Gamma$.
By construction of ${\cal X}_T$ and \eqref{stimaPF2}, we have that $\widehat s_{|_\Gamma},\widehat c_{|_\Gamma}\in  L^2(0,T;H^{3/2}(\Gamma))$.
Here, for the sake of simplicity, we use the same notation $\widehat c$ and $\widehat s$ on $\Gamma$ also for their traces.
Owing to standard theory for ordinary differential equations, exploiting \eqref{A7}
it is straightforward to find a unique $\widehat r={\cal S}_2(\widehat s,\widehat c)$ solving the Cauchy problem, for almost any $x \in \Gamma$,
\begin{equation}\label{eqPF2}
\partial_t r+\xi+\Psi'(r)+G(r,\widehat c,\widehat s)=F, \quad\hbox{a.e. in }(0,T), \quad  \quad r(0) = r_0.
\end{equation}
Testing  \eqref{eqPF2} (replacing $r$ with $\widehat r$)  by $\partial_t \widehat r$ and integrating over $(0,t)$, by means of Young's inequality,  the
smoothness conditions on $\Psi$ (see \eqref{A7}) and the chain rule for the sub-differential $\partial W$, we obtain (see \eqref{A5})
\begin{align}\label{stimarI}
&\frac 1 2\int_0^t\|\partial_t \widehat r\|^2_{L^2(\Gamma)}+\int_\Gamma W(\widehat r(t))\leq \int_\Gamma W(\widehat r(0))\\\no
&+ C\left(1+\int_0^t(\|\widehat s\|^2_{L^2(\Gamma)}+\|\widehat c\|^2_{L^2(\Gamma)})+\|F\|_{L^2(0,T;L^2(\Gamma))}
+\int_0^t\|\partial_t\widehat r\|^2_{L^2(0,s;L^2(\Gamma)}\right)\\\no
&\leq C\left(1+\int_0^t\left(\|\widehat s\|^2_{H^2(\Omega)}+\|\widehat c\|^2_{H^2(\Omega)}\right)+\int_0^t\|\partial_t\widehat r\|^2_{L^2(0,s;L^2(\Gamma)}\right).
\end{align}
Thus, using the Gronwall lemma and recalling \eqref{stimaPF1bis},   we eventually get
\begin{equation}\label{stimaPF3}
\|\widehat r\|_{H^1(0,T;L^2(\Gamma))}\leq C(R).
\end{equation}

\subsection{Third step: finding $s$}

We now consider the initial and boundary value problem \eqref{eq1w}, \eqref{eq1wbis}, \eqref{eq4w}$_1$, where $\widehat c={\cal S}_1(\widehat s)$
and $\widehat r={\cal S}_2(\widehat s,{\cal S}_1(\widehat s))$.
Applying standard results for linear parabolic equations,  we have that there exists a unique solution
$s={\cal S}_3(\widehat c,\widehat r)={\cal S}(\widehat s)$ belonging to $ H^1(0,T;H^1(\Omega))\cap L^2(0,T;H^2(\Omega))$.

We first test the resulting equation by $s$ and integrate over $(0,t)$. This gives
\begin{equation}\label{eqPFs1}
\int_0^t\int_\Omega\partial_t(\phi(\widehat c)s)s+\int_0^t\int_\Omega\phi(\widehat c)|\nabla s|^2+\int_0^t\int_\Gamma\nu(\widehat r)(s-\bar s)s
+\int_0^t\int_\Omega\phi(\widehat c)\widehat c s^2=0.
\end{equation}
We first deal with the first integral on the left hand side, integrating by parts in time and exploiting \eqref{sceltafi} and \eqref{eqPF1}. So that we obtain
\begin{align}\label{eqPFs2}
&\int_0^t\int_\Omega\partial_t(\phi(\widehat c)s)s=\int_\Omega\phi(\widehat c(t))s^2(t)-\int_\Omega\phi(c_0)s_0^2-\int_0^t\int_\Omega\phi(\widehat c)ss_t\\\no
&=\int_\Omega\phi(\widehat c(t))s^2(t)-\int_\Omega\phi(c_0)s_0^2-\frac 1 2\int_\Omega\phi(\widehat c(t))s^2(t) \\\no
&+\frac 1 2\int_\Omega\phi(c_0)s_0^2+\frac 1 2\int_0^t\int_\Omega B\widehat c_ts^2.
\end{align}
Combining \eqref{eqPFs1} with \eqref{eqPFs2}, we get (see \eqref{stimaPF0} and \eqref{stimaPF2})
\begin{align}\label{stimaPFs1}
&m\left( \frac{1}{2}\|s(t)\|^2_{L^2(\Omega)} + \int_0^t\|\nabla s\|^2_{L^2(\Omega)}\right) + \int_0^t\int_\Gamma \nu(\widehat r) s^2 \\ \no
&\leq \frac 12 M \|s_0\|^2_{L^2(\Omega)} + \frac12 \int_0^t\int_\Omega |B||\partial_t\widehat c|s^2
+ \int_0^t\int_\Omega \Phi(\widehat c) \widehat cs^2 + \int_0^t\int_\Gamma\nu(\widehat r)|s\bar s|\\\no
&\leq \frac 12 M \|s_0\|^2_{L^2(\Omega)} + \frac12 \int_0^t\int_\Gamma\nu(\widehat r)\bar s^2 \\\no
&+ \frac12 \int_0^t\int_\Omega (|B||\partial_t\widehat c|+ 2MC_0 )s^2 + \frac 12\int_0^t\int_\Gamma\nu(\widehat r)s^2,
\end{align}
from which we deduce
\begin{align}\label{stimaPFs1bis}
&m\left(\|s(t)\|^2_{L^2(\Omega)} + 2\int_0^t\|\nabla s(\sigma)\|^2_{L^2(\Omega)}\right) + \int_0^t\int_\Gamma \nu(\widehat r(\sigma)) s^2(\sigma) \\ \no
&\leq  M \|s_0\|^2_{L^2(\Omega)} +  \|\nu\|_{L^\infty(\mathbb{R})} \|\bar s\|^2_{L^\infty(0,T;L^2(\Gamma))} \\ \no
&+  \int_0^t(|B|\|\partial_t\widehat c(\sigma)\|_{L^\infty(\Omega)} + 2MC_0 )\|s(\sigma)\|^2_{L^2(\Omega)},
\end{align}
and then
\begin{align}\label{stimaPFs1ter}
&\|s(t)\|^2_{L^2(\Omega)}
\leq C +  \int_0^t C( \|\partial_t\widehat c(\sigma)\|_{L^\infty(\Omega)} + 1)\|s(\sigma)\|^2_{L^2(\Omega)}.
\end{align}
We can now use a generalized Gronwall's lemma (see \cite[Teorema 2.1]{baiocchi}).  Note that $\|\partial_t \widehat c\|_{L^\infty(\Omega)}$
is bounded in $L^2(0,T)$ by \eqref{stimaPF2}.  We obtain
\begin{equation}\label{stimaPFs2}
\|s\|_{L^\infty(0,T;L^2(\Omega))}\leq C(R),
\end{equation}
and, thanks to \eqref{stimaPFs1bis}, we get
\begin{equation}\label{stimaPFs2bis}
\|s\|_{L^\infty(0,T;L^2(\Omega))\cap L^2(0,T;H^1(\Omega))}\leq C(R).
\end{equation}
Recalling that $s$ is a strong solution, we test equation \eqref{eq1w}  by $\partial_t s$, where
$\widehat c={\cal S}_1(\widehat s)$ and
$\widehat r={\cal S}_2(\widehat s,{\cal S}_1(\widehat s))$.
Integrating over $(0,t)$ we obtain
\begin{align}\label{eqPFs3}
&\int_0^t\int_\Omega\phi(\widehat c)|\partial_t s|^2 +\int_0^t\int_\Omega \phi'(\widehat c)\partial_t \widehat c s\partial_t s
+\int_0^t\int_\Omega\phi(\widehat c)\nabla s\nabla \partial_ts\\\no
&+\int_0^t\int_\Gamma\nu(\widehat r)(s-\bar s)\partial_t s+\int_0^t\int_\Omega\phi(\widehat c) \widehat c s\partial_ts = 0.
\end{align}
Let us consider the identities
\begin{align}\label{eqI2s}
\int_0^t\int_\Omega\phi(\widehat c)\nabla s\nabla \partial_ts
&=\frac 1 2\int_\Omega\phi(\widehat c(t))|\nabla s(t)|^2-\frac 1 2\int_\Omega \phi(c_0)|\nabla s_0|^2 \\\no
&- \frac12\int_0^t\int_\Omega B\partial_t \widehat c|\nabla s|^2,
\end{align}
and
\begin{align}\label{eqI3s}
&\int_0^t\int_\Gamma\nu(\widehat r)(s-\bar s)\partial_t s\\\no
&=\frac 12\int_\Gamma\nu( \widehat r(t))(s-\bar s)^2(t)-\frac 12 \int_\Gamma \nu(r_0)(s_0-\bar s)^2
-\int_0^t\int_\Gamma\nu'(\widehat r)\partial_t \widehat r(s-\bar s)^2.
\end{align}
Then, on account of \eqref{eqI2s} and \eqref{eqI3s}, from \eqref{eqPFs3} we deduce
\begin{align}\label{eqPFs3bis}
&\int_0^t\int_\Omega\phi(\widehat c)|\partial_t s|^2 + \frac 1 2\int_\Omega\phi(\widehat c(t))|\nabla s(t)|^2 + \frac 12\int_\Gamma\nu(\widehat r(t))s^2(t)\\\no
&=\frac 12\int_0^t\int_\Omega B\partial_t \widehat c|\nabla s|^2- \int_0^t\int_\Omega \phi'(\widehat c)\partial_t \widehat c s\partial_t s
- \int_0^t\int_\Omega\phi(\widehat c) \widehat c s\partial_ts \\\no
&\frac 1 2\int_\Omega \phi(c_0)|\nabla s_0|^2 +\frac 12 \int_\Gamma \nu(r_0)(s_0-\bar s)^2
+\frac 12\int_0^t\int_\Gamma\nu'(\widehat r)\partial_t \widehat r(s-\bar s)^2 \\\no
&+ \int_\Gamma\nu(\widehat r(t))s\bar s(t) - \frac 12\int_\Gamma\nu(\widehat r(t))\bar s^2(t).
\end{align}
Observe now that
\begin{align}\label{eqI1s}
& \left |\int_0^t\int_\Omega \phi'(\widehat c)\partial_t \widehat c s\partial_t s \right |
\leq |B|\int_0^t\|\partial_t \widehat c\|_{L^\infty(\Omega)}\|s\|_{L^2(\Omega)}|\|\partial_ts\|_{L^2(\Omega)}\\\no
&\leq \frac m 2 \int_0^t\|\partial_t s\|^2_{L^2(\Omega)}+ \frac{B^2M^2C_0^2}{m}\int_0^t\|\widehat s\|^2_{L^\infty(\Omega)}\|s\|^2_{L^2(\Omega)}\\\no
&\leq  \frac m 2 \int_0^t\|\partial_t s\|^2_{L^2(\Omega)}+C(R),
\end{align}
due to the fact that $\|s\|^2_{L^2(\Omega)}$ is bounded in $L^\infty(0,T)$ (see \eqref{stimaPFs2}) and $\|\widehat s\|^2_{L^\infty(\Omega)}$ is bounded
in $L^1(0,T)$ by definition of ${\cal X}_{R,T}$.
Moreover, since  $\|\partial_t \widehat c\|_{L^\infty(\Omega)}$ is bounded in $L^2(0,T)$ (see \eqref{stimaPF2}), it holds (cf. \eqref{stimaPFs2bis})
\begin{equation}
\left |\int_0^t\int_\Omega B\partial_t \widehat c|\nabla s|^2 \right |\leq |B|\int_0^t\|\partial_t \widehat c\|_{L^\infty(\Omega)}\|\nabla s\|^2_{L^2(\Omega)}
\leq C(R).
\end{equation}
Recalling  \eqref{stimaPF3} and \eqref{stimaPFs2bis} we get
\begin{equation}
\left |\int_0^t\int_\Gamma\nu'(\widehat r)\partial_t \widehat r(s-\bar s)^2\right |
\leq C\int_0^t\|\partial_t \widehat r\|_{L^2(\Gamma)}\|s-\bar s\|^2_{L^4(\Gamma)} \leq C(R).
\end{equation}
Analogously, we have (see \eqref{stimaPFs2bis})
\begin{equation}
\left |\int_0^t\int_\Gamma\nu(\widehat r)s\bar s\right |\leq C\int_0^t\left ( \|\bar s\|^2_{L^2(\Gamma)}+\|s\|^2_{H^1(\Omega)}\right ) \leq C + C(R).
\end{equation}
Thus we deduce (see \eqref{stimaPFs2})
\begin{align}\label{eqI4s}
&\int_0^t\int_\Omega\phi(\widehat c)\widehat cs\partial_ts\leq MC_0\int_0^t\|s\|_{L^2(\Omega)}\|\partial_ts\|_{L^2(\Omega)}\\\no
&\leq \frac m4\int_0^t\|\partial_ts\|^2_{L^2(\Omega)}+\frac{M^2C_0^2}{m}\int_0^t\|s\|^2_{L^2(\Omega)} \leq \frac m4\int_0^t\|\partial_ts\|^2_{L^2(\Omega)}+C(R).
\end{align}
Then, combining \eqref{eqI1s}--\eqref{eqI4s} with \eqref{eqPFs3bis}, we find
\begin{align}\label{eqPFs3ter}
&\frac m4\|\partial_t s\|_{L^2(0,T;L^2(\Omega))}^2 + \frac m2\|\nabla s\|_{L^\infty(0,T;L^2(\Omega))}^2  \leq C +  C(R),
\end{align}
which gives (see also  \eqref{stimaPFs2})
\begin{equation}\label{stimaPFs3}
\|\partial_t s\|_{L^2(0,T;L^2(\Omega))}+\|s\|_{L^\infty(0,T;H^1(\Omega))}\leq C(R).
\end{equation}
Actually, we can now estimate $\partial_t s$ in $L^\infty(0,T;L^2(\Omega))$. To this aim let us (formally)  differentiate the equation for $s$
with respect to time and get
\begin{align}\label{eqsder}
&\phi(\widehat c) s_{tt}+2B\widehat c_t s_t+B\widehat c_{tt}s-\hbox{div }(\phi(\widehat c)\nabla s)_t\\\no
&=-B\widehat c_t\widehat cs-\phi(\widehat c)\widehat c_ts-\phi(\widehat c)\widehat cs_t.
\end{align}
Observe that we can determine  $s_1=s_t(0) $ from the identity (cf. \eqref{reazdiff})
\begin{equation}-B\phi(c_0)s^2_0c_0 + \phi(c_0)s_1 - {\rm div}(\phi(c_0)\nabla s_0) = \phi(c_0)s_0c_0,\end{equation}
and, due to the regularity of $s_0$ and $c_0$, we can deduce that $s_1 \in L^2(\Omega)$.
The boundary condition turns out to be rewritten as (see Remark \ref{stazs})
\begin{equation}\label{boundsder}
(\phi(\widehat c)\partial_n s)_t = \phi(\widehat c)\partial_ns_t+B\widehat c_t\partial_ns= - \nu'(\widehat r)\widehat r_t(s-\bar s) - \nu(\widehat r)s_t.
\end{equation}
In order to make our argument rigorous, we can observe, for instance, that the resulting problem \eqref{eqsder}--\eqref{boundsder} (written with respect to
the unknown $\tilde s=s_t$) has the  same structure of our original problem, and thus it admits a unique solution in the standard energy space.

We thus test \eqref{eqsder} by $s_t$ and integrate over $(0,t)$, obtaining
\begin{align}\label{eqsder2}
&\frac12 \int_\Omega\phi(\widehat c) s^2_{t}(t) - \frac12 \int_\Omega\phi(c_0) s^2_1
-\frac12 \int_0^t\int_\Omega B \widehat c_t s^2_{t}
+2B\int_0^t\int_\Omega\widehat c_t s^2_t\\\no
&+ B\int_0^t\int_\Omega\widehat c_{tt}ss_t -\int_0^t\int_\Gamma(\phi(\widehat c)\partial_n s)_ts_t +
\int_0^t\int_\Omega(\phi(\widehat c)\nabla s)_t\cdot \nabla s_t \\\no
&=-B\int_0^t\int_\Omega\widehat c_t\widehat css_t-\int_0^t\int_\Omega\phi(\widehat c)\widehat c_tss_t-\int_0^t\int_\Omega\phi(\widehat c)\widehat cs^2_t,
\end{align}
from which, recalling \eqref{eqPF1} and \eqref{boundsder}, we deduce the identity
\begin{align}\label{eqsder3}
&\frac12 \int_\Omega\phi(\widehat c) s^2_{t}(t)
+\frac32 \int_0^t\int_\Omega B \widehat c_t s^2_{t}
- B\int_0^t\int_\Omega ( \phi(\widehat c)s\widehat c )_{t}ss_t + \int_0^t\int_\Gamma\nu(\widehat r)s^2_t\\\no
&+\int_0^t\int_\Gamma\nu'(\widehat r)\widehat r_t(s-\bar s)s_t
+ \int_0^t\int_\Omega \phi(\widehat c)|\nabla s_t|^2  + \int_0^t\int_\Omega B \widehat c_t\nabla s \cdot \nabla s_t \\\no
&=\frac12 \int_\Omega\phi(c_0) s^2_{1}-B\int_0^t\int_\Omega\widehat c_t\widehat css_t-\int_0^t\int_\Omega\phi(\widehat c)\widehat c_tss_t
-\int_0^t\int_\Omega\phi(\widehat c)\widehat cs^2_t.
\end{align}
Recalling the properties of $\widehat c$, $\widehat r$ and $s$, from \eqref{eqsder3} we obtain the estimate
\begin{align}\label{stimasder}
&\frac 12\int_\Omega\phi(\widehat c(t))|s_t(t)|^2+\int_0^t\int_\Omega\phi(\widehat c)|\nabla s_t|^2+\int_0^t\int_\Gamma\nu(\widehat r)|s_t|^2\\\no
&\leq \frac 12\int_\Omega\phi(c_0)s_1^2+|B| \int_0^t\|\widehat c_t\|_{L^\infty(\Omega)}\|\nabla s\|_{L^2(\Omega)}\|\nabla s_t\|_{L^2(\Omega)}\\\no
&+\int_0^t\|\nu'\|_{L^\infty(\mathbb{R})}\|\widehat r_t\|_{L^2(\Gamma)}\|s-\bar s\|_{L^4(\Gamma)}\|s_t\|_{L^4(\Gamma)}\\\no
&+ \frac32|B|\int_0^t\|\widehat c_t\|_{L^\infty(\Omega)}\|s_t\|^2_{L^2(\Omega)}
+|B|\int_0^t\|\widehat c_t\|_{L^\infty(\Omega)}\|\widehat c\|_{L^\infty(\Omega)}\|s\|_{L^2(\Omega)}\|s_t\|_{L^2(\Omega)}\\\no
&+M\int_0^t\|\widehat c_t\|_{L^\infty(\Omega)}\|s\|_{L^2(\Omega)}\|s_t\|_{L^2(\Omega)}
+MC_0\int_0^t\|s_t\|^2_{L^2(\Omega)} + B\int_0^t\int_\Omega ( \phi(\widehat c)s\widehat c )_{t}ss_t.
\end{align}
On account of the Young inequality and the Sobolev injection $H^1(\Omega)\hookrightarrow L^4(\Omega)$, observe that it holds
\begin{align}\label{stimasder2}
& B\int_0^t\int_\Omega ( \phi(\widehat c)s\widehat c )_{t}ss_t = B\int_0^t\int_\Omega \left ( B \widehat c_t \widehat c s^2s_t
+  \phi(\widehat c)\widehat c ss^2_t +  \phi(\widehat c)\widehat c_{t}s^2s_t \right ) \\\no
&\leq C\int_0^t \left( \| \widehat c_t\|_{L^\infty(\Omega)}\| s\|^2_{L^4(\Omega)}\| s_t\|_{L^2(\Omega)}
+ \| s\|_{L^4(\Omega)}\| s_t\|_{L^4(\Omega)}\| s_t\|_{L^2(\Omega)} \right)\\\no
&\leq C\int_0^t \left( \| \widehat c_t\|_{L^\infty(\Omega)}\| s\|^2_{L^4(\Omega)}\| s_t\|_{L^2(\Omega)}
+ \| s\|_{L^4(\Omega)}\| s_t\|_{H^1(\Omega)}\| s_t\|_{L^2(\Omega)} \right)\\\no
&\leq C\int_0^t \left( \| \widehat c_t\|_{L^\infty(\Omega)}\| s\|^2_{L^4(\Omega)}\| s_t\|_{L^2(\Omega)}
+ \| s\|_{L^4(\Omega)}\| s_t\|_{L^2(\Omega)}\| s_t\|_{L^2(\Omega)} \right)\\\no
& + C\int_0^t \| s\|_{L^4(\Omega)} \| \nabla s_t\|_{L^2(\Omega)} \| s_t\|_{L^2(\Omega)} \\\no
&\leq C\int_0^t \left( \| \widehat c_t\|_{L^\infty(\Omega)}\| s\|^2_{L^4(\Omega)}\| s_t\|_{L^2(\Omega)}
+ \| s\|_{L^4(\Omega)} \| s_t\|^2_{L^2(\Omega)} + \| s\|^2_{L^4(\Omega)} \| s_t\|^2_{L^2(\Omega)} \right)\\\no
&+\frac m4\int_0^t\| \nabla s_t\|^2_{L^2(\Omega)}.
\end{align}
Combining \eqref{stimasder} and \eqref{stimasder2} we then find
\begin{align}\label{stimasder3}
&\frac m2\|s_t(t)\|_{L^2(\Omega)}^2+m\int_0^t\|\nabla s_t\|^2_{L^2(\Omega)}\\\no
&\leq C + C \int_0^t \Big ( \|\widehat c_t\|^2_{L^\infty(\Omega)}\|\nabla s\|^2_{L^2(\Omega)} \Big)  +  \frac m8 \int_0^t \|\nabla s_t\|^2_{L^2(\Omega)}\\\no
&+ C \int_0^t \Big (\|\widehat r_t\|^2_{L^2(\Gamma)}(\|s\|^2_{L^4(\Gamma)}+\|\bar s\|^2_{L^4(\Gamma)}) + \|s_t\|^2_{L^2(\Omega)} \Big)
+ \frac m8 \int_0^t  \|\nabla s_t\|^2_{L^2(\Omega)}\\\no
&+C\int_0^t\Big (\|\widehat c_t\|_{L^\infty(\Omega)}\|s_t\|^2_{L^2(\Omega)}
+\|\widehat c_t\|_{L^\infty(\Omega)}\|s\|_{L^2(\Omega)}\|s_t\|_{L^2(\Omega)} + \|s\|^2_{L^2(\Omega)}\Big )\\\no
&+C\int_0^t\Big (\|\widehat c_t\|_{L^\infty(\Omega)}\|s\|^2_{L^4(\Omega)}\|s_t\|_{L^2(\Omega)}
+ \|s\|^2_{L^4(\Omega)}\|s_t\|^2_{L^2(\Omega)}+ \|s\|_{L^4(\Omega)}\|s_t\|^2_{L^2(\Omega)}\Big )\\\no
& + \frac m4  \int_0^t \|\nabla s_t\|^2_{L^2(\Omega)}.
\end{align}
Here, we have used the Young inequality, the trace theorem and the Sobolev injection once more.
In addition note that, by \eqref{stimaPF2}, \eqref{stimaPF3},
and \eqref{stimaPFs3}, from \eqref{stimasder3}  we deduce
\begin{align}\label{stimasder4}
&\frac m2\|s_t(t)\|_{L^2(\Omega)}^2 + \frac m2  \int_0^t \|\nabla s_t\|^2_{L^2(\Omega)} \\\no
&\leq C + C(R) +  C\int_0^t \Big ( 1 + \|\widehat c_t\|_{L^\infty(\Omega)} + \|s\|^2_{L^4(\Omega)}+ \| s\|_{L^4(\Omega)} \Big)\|s_t\|^2_{L^2(\Omega)}\\\no
& + C\int_0^t \|\widehat c_t\|_{L^\infty(\Omega)} \big ( \|s\|_{L^2(\Omega)}+ \| s\|^2_{L^4(\Omega)} \big)\|s_t\|_{L^2(\Omega)}.
\end{align}
Then we can apply the generalized Gronwall lemma, getting
\begin{equation}
\|s_t\|_{L^\infty(0,T;L^2(\Omega))}\leq C(R),
\end{equation}
which, along with \eqref{stimasder4}, gives
\begin{equation}\label{boundst}
\|s_t\|_{L^\infty(0,T;L^2(\Omega))\cap L^2(0,T;H^1(\Omega))}\leq C(R).
\end{equation}
Let us now rewrite  equation \eqref{eq1w} in the following form
\begin{equation}\label{eq1strong}
-\hbox{div }(\phi(\hat c)\nabla s)=-\phi(\hat c) \hat c s-(\phi(\hat c)s)_t,
\end{equation}
and combine it with the boundary condition on $\Gamma\times(0,T)$
\begin{equation}\label{eqbordoforte}
\partial_ns=-\frac 1{\phi(\hat c)}\nu(\hat r)(s-\bar s).
\end{equation}
We first observe that the right-hand side in \eqref{eq1strong} is bounded in $L^\infty(0,T;L^2(\Omega))$ and the right-hand side of \eqref{eqbordoforte}
is bounded in $L^\infty(0,T;H^{1/2}(\Gamma))$ (see \eqref{A3} and \eqref{stimaPFs3}). Thus, using standard ellliptic regularity results
(see, e.g., \cite[Chap. 5]{Taylor}), we get
\begin{equation}\label{stimaPFs4}
\|s\|_{L^\infty(0,T;H^2(\Omega))}\leq C(R).
\end{equation}
In order to prove that ${\cal S}:{\cal X}_{R,T}\rightarrow{\cal X}_{R,T}$, we have to show that $s\geq0$ a.e. in $(0,T)\times \Omega$ and
$\|s\|_{L^2(0,T;H^2(\Omega))}\leq R$.
We deduce the positivity of $s$ by applying a maximum principle argument (see, for instance, \cite[Lemma 2.1]{pao}). Indeed $s$ solves the problem
\begin{align}
 &\phi(\hat c) \partial_t s-\hbox{div }(\phi(\hat c)\nabla s) + \phi(\hat c) \hat c s = B\phi(\hat c)\hat c s^2 \geq 0,\quad\hbox{in }\Omega, \\
&\phi(c)\partial_ns+ \nu(r)s_{|_\Gamma} = \nu(r)\bar s \geq 0,\quad\hbox{on }\Gamma,\\
&s(0) = s_0 \geq 0,\quad\hbox{in }\Omega.
\end{align}
So that we have
\begin{equation}\label{spos}
s(x,t)\geq 0,\quad\hbox{for a.e. $(x,t)\in \Omega\times(0,T)$}.
\end{equation}
Then, let us exploit \eqref{boundst} and \eqref{stimaPFs4}  to get
\begin{align}
&\|s\|_{L^\infty(0,T;L^2(\Omega))\cap L^2(0,T;H^2(\Omega))}\\\no
&\leq T\|s_t\|_{L^\infty(0,T;L^2(\Omega))}+\|s_0\|_{L^2(\Omega)}+T^{1/2}\|s\|_{L^\infty(0,T;H^2(\Omega))}\\\no
&\leq (T+T^{1/2})C (R)+\|s_0\|_{L^2(\Omega)}.
\end{align}
Assuming $\widehat T \in (0,T]$ sufficiently small such that
\begin{equation}
 (\widehat T+\widehat T^{1/2})C (R)\leq R-\|s_0\|_{L^2(\Omega)},
 \end{equation}
we infer
\begin{equation}
\|s\|_{L^\infty(0,T;L^2(\Omega))\cap L^2(0,T;H^2(\Omega))}\leq R,
\end{equation}
so that  ${\cal S}:{\cal X}_{R,\widehat T}\rightarrow{\cal X}_{R,\widehat T}$. It is clear from the construction that
a fixed point of ${\cal S}$ is a local (in time) solution to our problem.

\subsection{The contraction argument}

Let us fix $\widehat s_i\in{\cal X}_{R,\widehat T}$, $i=1,2$, and set $\widehat c_i={\cal S}_1(\widehat s_i)$,
$\widehat r_i={\cal S}_2(\widehat s_i, \widehat c_i)$ and $s_i={\cal S}_3(\widehat c_i, \widehat r_i) = {\cal S}(\widehat s_i)$.
We first estimate $(\widehat c_1-\widehat c_2)$ and $(\widehat r_1-\widehat r_2)$ in suitable norms. To this aim, let us first write \eqref{eqPF1} for $i=1,2$,
take the difference and test by $\widehat c_1-\widehat c_2$. After integrating over $(0,t)$  and using the Young inequality, we get (see \eqref{stimaPF1})
{\ele \begin{align}\label{contr1}
&\frac 1 2\int_\Omega|(\widehat c_1-\widehat c_2)(t)|^2 \\\no
&= - \int_0^t\int_\Omega \phi(\widehat c_1)\widehat c_1(\widehat s_1-\widehat s_2)(\widehat c_1-\widehat c_2)
-\int_0^t\int_\Omega(\phi(\widehat c_1)\widehat c_1-\phi(\widehat c_2)\widehat c_2)\widehat s_2(\widehat c_1-\widehat c_2)\\\no
&=-\int_0^t\int_\Omega \phi(\widehat c_1)\widehat c_1(\widehat s_1-\widehat s_2)(\widehat c_1-\widehat c_2)
-\int_0^t\int_\Omega\big ( A + B(\widehat c_1+\widehat c_2) \big )\widehat s_2(\widehat c_1-\widehat c_2)^2\\\no
&\leq C\left(\int_0^t\|\widehat s_1-\widehat s_2\|^2_{L^2(\Omega)}+\|\widehat c_1-\widehat c_2\|^2_{L^2(\Omega)}
+\int_0^t\|\widehat s_2\|_{L^\infty(\Omega)}\|\widehat c_1-\widehat c_2\|^2_{L^2(\Omega)}\right).
\end{align}}
Using once more the generalized Gronwall lemma and the fact that $\|\widehat s_2\|_{L^\infty(\Omega)}$ is bounded in $L^2(0,\widehat T)$
(see \eqref{stimaPFs4}), we obtain
\begin{equation}\label{contr2}
\|(\widehat c_1-\widehat c_2)(t)\|_{L^2(\Omega)}\leq C(R)\|\widehat s_1-\widehat s_2\|_{L^2(0,t;L^2(\Omega))}.
\end{equation}
As a consequence, we deduce
\begin{align}\label{contr3}
&\|(\widehat c_1-\widehat c_2)_t\|_{L^2(0,t;L^2(\Omega))}\leq \|\phi(\widehat c_1)\widehat c_1(\widehat s_1-\widehat s_2)\|_{L^2(0,t;L^2(\Omega))}\\\no
&+\|\widehat s_2(A(\widehat c_1-\widehat c_2)+B(\widehat c_1+\widehat c_2)(\widehat c_1-\widehat c_2))\|_{L^2(0,t;L^2(\Omega))}\\\no
&\leq C\left(\|\widehat s_1-\widehat s_2\|_{L^2(0,t;L^2(\Omega))}
+\|\widehat s_2\|_{L^2(0,t;L^\infty(\Omega))}\|\widehat c_1-\widehat c_2\|_{L^\infty(0,t;L^2(\Omega))}\right)\\\no
&\leq C(R)\|\widehat s_1-\widehat s_2\|_{L^2(0,t;L^2(\Omega))}.
\end{align}
Let us estimate $\|\widehat c_1-\widehat c_2\|_{L^2(0,\widehat T;H^1(\Omega))}$. To this aim, we proceed by
integrating  the  difference of equation \eqref{eqPF1} written for the two indices. This yields
\begin{align}\label{contr4}
&\|(\widehat c_1-\widehat c_2)(t)\|_{H^1(\Omega)}
\leq\int_0^t\|\phi(\widehat c_1)\widehat c_1\|_{L^\infty(\Omega)}\|\widehat s_1-\widehat s_2\|_{H^1(\Omega)}\\\no
&+\int_0^t(|B|\|\widehat c_1\|_{L^\infty(\Omega)}
+ \|\phi(\widehat c_1)\|_{L^\infty(\Omega)})\|\nabla \widehat c_1\|_{L^2(\Omega)}\|\widehat s_1-\widehat s_2\|_{L^2(\Omega)}\\\no
&+\int_0^t\|\phi(\widehat c_1)\widehat c_1 - \phi(\widehat c_2)\widehat c_2\|_{L^\infty(\Omega)}\|\widehat s_2\|_{L^2(\Omega)}
\|\widehat c_1-\widehat c_2\|_{H^1(\Omega)}\\\no
&+\int_0^t\|B\widehat c_1\nabla \widehat c_1 +\phi(\widehat c_1)\nabla \widehat c_1 - B\widehat c_2\nabla \widehat c_2
- \phi(\widehat c_2)\nabla \widehat c_2\|_{L^2(\Omega)}\|\widehat s_2\|_{L^\infty(\Omega)}\|\widehat c_1-\widehat c_2\|_{L^2(\Omega)}\\\no
&\leq C \int_0^t\|\widehat s_1-\widehat s_2\|_{H^1(\Omega)} + C \int_0^t\|\nabla \widehat c_1\|_{L^2(\Omega)}\|\widehat s_1-\widehat s_2\|_{L^2(\Omega)}\\\no
&+C \int_0^t\|s_2\|_{L^2(\Omega)} \|\widehat c_1-\widehat c_2\|_{H^1(\Omega)} \\\no
&+ C \int_0^t(\|\nabla \widehat c_1\|_{L^2(\Omega)}
+ \|\nabla \widehat c_2\|_{L^2(\Omega)})\|\widehat s_2\|_{L^\infty(\Omega)} \|\widehat c_1-\widehat c_2\|_{L^2(\Omega)}.
\end{align}
An application of the Gronwall lemma and the previous estimates give
\begin{equation}\label{contr5}
\|(\widehat c_1-\widehat c_2)(t)\|_{H^1(\Omega)}\leq C(R)\int_0^t\|\widehat s_1-\widehat s_2\|_{H^1(\Omega)}.
\end{equation}
In order to estimate $(\widehat r_1- \widehat r_2)$, we write  \eqref{eqPF2} for the two indices, take the difference and test by $(\widehat r_1-\widehat r_2)$.
Then, we integrate in time and exploit the monotonicity of $W'$ and \eqref{A7}. This implies
\begin{equation}
\frac 1 2\|(\widehat r_1-\widehat r_2)(t)\|^2_{L^2(\Gamma)}
\leq C\int_0^t\|\widehat c_1-\widehat c_2\|^2_{L^2(\Gamma)}+\|\widehat s_1-\widehat s_2\|^2_{L^2(\Gamma)}+\|\widehat r_1-\widehat r_2\|^2_{L^2(\Gamma)}.
\end{equation}
Thus, due to \eqref{contr5}, we can eventually deduce
\begin{equation}\label{contr6}
\|(\widehat r_1-\widehat r_2)(t)\|_{L^2(\Gamma)}\leq C(R) \|\widehat s_1-\widehat s_2\|_{L^2(0,t;H^1(\Omega))}.
\end{equation}
Consider now the two equations for $s_1$ and $s_2$. Subtracting the second from the first, we get
\begin{align}\label{contr8}
&(\phi(\widehat c_1)s_1)_t-(\phi(\widehat c_2)s_2)_t-\hbox{div }(\phi(\widehat c_1)\nabla s_1-\phi(\widehat c_2)\nabla s_2)\\\no
&=-\phi(\widehat c_1)\widehat c_1s_1+\phi(\widehat c_2)\widehat c_2s_2.
\end{align}
Recalling the expression for $\phi$, we obtain
\begin{align}\label{contr9}
&\phi(\widehat c_1)(s_1-s_2)_t+B(\widehat c_1-\widehat c_2)\partial_ts_2+B\partial_t\widehat c_1(s_1-s_2)+B(\widehat c_1-\widehat c_2)_ts_2\\\no
&-\hbox{div }(\phi(\widehat c_1)\nabla(s_1-s_2)+B(\widehat c_1-\widehat c_2)\nabla s_2)\\\no
&=-\phi(\widehat c_1)\widehat c_1(s_1-s_2)-(A(\widehat c_1-\widehat c_2)+B(\widehat c_1-\widehat c_2)(\widehat c_1+\widehat c_2))s_2,
\end{align}
with the boundary condition
\begin{equation}\label{contr10}
\phi(\widehat c_1)\partial_n s_1-\phi(\widehat c_2)\partial_ns_2=-\nu(\widehat r_1)(s_1-\bar s)+\nu(\widehat r_2)(s_2-\bar s).
\end{equation}
Testing \eqref{contr9} by $s_1-s_2$ and integrating over $(0,t)$, we find
\begin{align}\label{contr11}
&\frac 1 2\int_0^t\int_\Omega\phi(\widehat c_1)\frac{d}{dt}|s_1-s_2|^2+B\int_0^t\int_\Omega (\widehat c_1-\widehat c_2)\partial_ts_2(s_1-s_2)\\\no
&+B\int_0^t\int_\Omega\partial_t\widehat c_1(s_1-s_2)^2 + B\int_0^t\int_\Omega (\widehat c_1-\widehat c_2)_ts_2(s_1-s_2)\\\no
&+\int_0^t\int_\Omega \phi(\widehat c_1)|\nabla(s_1-s_2)|^2+B\int_0^t\int_\Omega  (\widehat c_1-\widehat c_2)\nabla s_2\nabla(s_1-s_2)\\\no
&+\int_0^t\int_\Gamma\nu(\widehat r_1)(s_1-s_2)^2+(\nu(\widehat r_1)-\nu(\widehat r_2))(s_2-\bar s)(s_1-s_2)\\\no
&=-\int_0^t\int_\Omega\phi(\widehat c_1)\widehat c_1(s_1-s_2)^2\\\no
&-\int_0^t\int_\Omega \big (A(\widehat c_1-\widehat c_2)+B(\widehat c_1-\widehat c_2)(\widehat c_1+\widehat c_2)\big )s_2(s_1-s_2).
\end{align}
Observe that
\begin{equation}\label{contr12}
\frac 1 2\int_0^t\int_\Omega\phi(\widehat c_1)\frac{d}{dt}|s_1-s_2|^2
=\frac 12\int_\Omega\phi(\widehat c_1(t))|(s_1-s_2)(t)|^2-\frac 12\int_0^t\int_\Omega B\partial_t \widehat c_1|s_1-s_2|^2,
\end{equation}
where
\begin{equation}\label{contr13}
\left |-\frac 12\int_0^t\int_\Omega B\partial_t \widehat c_1|s_1-s_2|^2 \right |
\leq C\int_0^t\|\partial_t\widehat c_1\|_{L^\infty(\Omega)}\|s_1-s_2\|^2_{L^2(\Omega)},
\end{equation}
and $\|\partial_t\widehat c_1\|_{L^\infty(\Omega)}$ is bounded in $L^2(0,\widehat T)$ (see \eqref{stimaPF2}).
Then, it holds (see \eqref{stimaPF2} and \eqref{stimaPFs3})
\begin{align}\label{contr14}
&\left | B\int_0^t\int_\Omega (\widehat c_1-\widehat c_2)\partial_ts_2(s_1-s_2) \right |
\leq C\int_0^t\|\partial_ts_2\|_{L^2(\Omega)}\|\widehat c_1-\widehat c_2\|_{L^4(\Omega)}\|s_1-s_2\|_{L^4(\Omega)}\\\no
&\leq \delta\int_0^t\|s_1-s_2\|^2_{H^1(\Omega)} +C\int_0^t\|\partial_t s_2\|^2_{L^2(\Omega)}\|\widehat c_1-\widehat c_2\|^2_{H^1(\Omega)}\\\no
&\leq \delta\int_0^t\|s_1-s_2\|^2_{H^1(\Omega)} +C(R)\int_0^t\|\partial_t s_2\|^2_{L^2(\Omega)} \|\widehat s_1-\widehat s_2\|^2_{L(0,s;H^1(\Omega))},
\end{align}
for a sufficiently small $\delta$ to be chosen later. Analogously, we have (see \eqref{contr3})
\begin{align}\label{contr15}
&\left | B\int_0^t\int_\Omega (\widehat c_1-\widehat c_2)_ts_2(s_1-s_2) \right |
\leq C\int_0^t\|s_2\|_{L^\infty(\Omega)}\|(\widehat c_1-\widehat c_2)_t\|_{L^2(\Omega)}\|s_1-s_2\|_{L^2(\Omega)}\\\no
&\leq \delta\|s_1-s_2\|^2_{L^\infty(0,t;L^2(\Omega))}+C\int_0^t\|(\widehat c_1-\widehat c_2)_t\|^2_{L^2(\Omega)}\\\no
&\leq \delta\|s_1-s_2\|^2_{L^\infty(0,t;L^2(\Omega))}+C(R)\|(\widehat s_1-\widehat s_2)\|^2_{L^2(0,t;L^2(\Omega))}.
\end{align}
The gradient terms are estimated using \eqref{contr5}, namely,
\begin{align}\label{contr16}
&\left | B\int_0^t\int_\Omega  (\widehat c_1-\widehat c_2)\nabla s_2\nabla(s_1-s_2) \right |\\\no
&\leq B\int_0^t\|\widehat c_1-\widehat c_2\|_{L^4(\Omega)}\|\nabla s_2\|_{L^4(\Omega)}\|\nabla(s_1-s_2)\|_{L^2(\Omega)}\\\no
&\leq \delta\int_0^t\|s_1-s_2\|^2_{H^1(\Omega)}+C(R)\int_0^t\|\widehat s_1-\widehat s_2\|^2_{L^2(0,s;H^1(\Omega))}.
\end{align}
Observe that we can control the boundary term in this way
\begin{align}\label{contr17}
&\left | \int_0^t\int_\Gamma(\nu(\widehat r_1)-\nu(\widehat r_2))(s_2-\bar s)(s_1-s_2)\right |\\\no
&\leq \|\nu'\|_{L^\infty(\mathbb{R})}\int_0^t\|\widehat r_1- \widehat r_2\|_{L^2(\Gamma)}\|s_2-\bar s\|_{L^4(\Gamma)}\|s_1-s_2\|_{L^4(\Gamma)}\\\no
&\leq\delta \int_0^t\|s_1-s_2\|^2_{H^1(\Omega)}+C\int_0^t\|r_1-r_2\|^2_{L^2(\Gamma)}\|s_2-\bar s\|^2_{L^4(\Gamma)}\\\no
&\leq\delta \int_0^t\|s_1-s_2\|^2_{H^1(\Omega)}+C(R)\int_0^t\|\widehat s_1-\widehat s_2\|^2_{L^2(0,s;H^1(\Omega))}.
\end{align}
{\ele Combining \eqref{contr11}--\eqref{contr17} and taking $\delta$ sufficiently small, we deduce
\begin{align}
&\|(s_1-s_2)(t)\|^2_{L^2(\Omega)}+\int_0^t\|(s_1-s_2)\|^2_{H^1(\Omega)}\\\no
&\leq C(R)\int_0^t(1+\|\partial_t s_2\|^2_{L^2(\Omega)})\|\widehat s_1-\widehat s_2\|^2_{L^2(0,s;H^1(\Omega))}\\\no
&+C(R)\|(\widehat s_1-\widehat s_2)\|^2_{L^2(0,t;L^2(\Omega))}+C\int_0^t\|\partial_t\widehat c_1\|_{L^\infty(\Omega)}\| s_1-s_2\|^2_{L^2(\Omega)}.
\end{align}
Finally, making use of the generalized Gronwall lemma, we get
\begin{align}\label{contr18}
&\|(s_1-s_2)(t)\|^2_{L^2(\Omega)}+\int_0^t\|s_1-s_2\|^2_{H^1(\Omega)}\\\no
&\leq C(R)\left(\int_0^t\|\widehat s_1-\widehat s_2\|^2_{L^2(\Omega)}+\int_0^t\|\widehat s_1-\widehat s_2\|^2_{L^2(0,s;H^1(\Omega))}\right).
\end{align}
Thus, for a suitable large $j\in \mathbb{N}$ (we recall that $R$ is fixed)} we have that ${\cal S}^j$ is a contraction in $\mathcal{X}_{R,\widehat T}$. Hence,
$\mathcal{S}$ admits a unique fixed point $s$, which is the local in time solution to our original problem.

\section{Global existence result for the initial problem}

In this section, we show that \eqref{A9} allows us to obtain a global a priori $L^\infty$ bound on $s$. Therefore, the previous fixed point argument
can now be carried out without restriction on $T$ so that the unique fixed point is solution to our problem on the whole time interval $(0,T)$.

\subsection{Uniform a priori bound for $s$}

Assuming  \eqref{A9}, let us prove that $s$ is globally bounded in $L^\infty(\Omega\times (0,T))$.
Let us recall that \eqref{spos} holds. Then,
we introduce the auxiliary function
\begin{equation}\label{defz}
z(x,t)=S_0-s(x,t), \quad\hbox{ for a.e. } (x,t)\in\Omega\times(0,T),
\end{equation}
and prove that $z\geq 0$. We observe that $z$ solves the problem
\begin{align}\label{eqz1}
&\partial_t(\phi(c)z)-\hbox{div }(\phi(c)\nabla z)+\phi(c)cz=(\phi(c)c+\phi'(c)c_t)S_0,\quad\hbox{ a.e. in }\Omega\times (0,T),\\\label{cauchyz}
&z(x,0)=S_0-s_0(x),\quad\hbox{ in }\Omega,\\\label{boundz}
&\phi(c)\partial_n z + \nu(r)z_{|_{\Gamma}} = \nu (r)(S_0 - \bar s), \quad\hbox{ a.e. on }\Gamma\times (0,T).
\end{align}
Exploiting \eqref{eq2w} (and \eqref{sceltafi})  we rewrite \eqref{eqz1} as follows
\begin{equation}
\partial_t(\phi(c)z)-\hbox{div }(\phi(c)\nabla z)+\phi(c)cz(1-BS_0)=S_0\phi(c)c(1-BS_0).
\end{equation}
Thanks to \eqref{A9}, we observe that
$$S_0\phi(c)c(1-BS_0) \geq 0, \quad S_0-s_0 \geq 0, \quad \nu (r)(S_0 - \bar s) \geq 0.$$
Then we can apply once more \cite[Lemma 2.1]{pao} and deduce that $z\geq 0$ in $\Omega\times(0,T)$. So that (see \eqref{spos})
\begin{equation}\label{bounds}
0\leq s(x,t)\leq S_0,\quad\hbox{ for a.e. }(x,t)\in \Omega\times (0,T).
\end{equation}

\subsection{Global a priori estimates}

Here we show that using \eqref{bounds} and the uniform bound for $c=\widehat c$ (see \eqref{stimaPF1}), then $\|s\|_{L^\infty(0,T;H^2(\Omega))}$ is bounded
by a constant depending only on the data.

First, we observe that  \eqref{eq2w} and \eqref{bounds} lead to a uniform estimate on $\partial_t c$
\begin{equation}
\|\partial_t c\|_{L^\infty(\Omega\times(0,T))}\leq C.
\end{equation}
Here $C$ also depends on $S_0$ (see \eqref{A9}).
Then, we can test \eqref{eq1w} by $s$ and integrate over $(0,t)$.
Arguing as in \eqref{stimaPFs1ter}--\eqref{stimaPFs2} we now deduce, in particular, that
\begin{equation}
\|s\|_{L^2(0,T;H^1(\Omega))}\leq C,
\end{equation}
where now $C$ does not depend on $R$.

Integrating \eqref{eq2w} with respect to time, and taking then the spatial gradient, we obtain the identity
$$\nabla c(t) = \nabla c_0 - \int_0^t( Bcs \nabla c + \phi(c) s \nabla c + \phi(c)c \nabla s).$$
By means of the previous estimates and the generalized Gronwall lemma we also get that
\begin{equation}
\|c\|_{L^\infty(0,T;H^1(\Omega))}\leq C,
\end{equation}
where $C$ does not depend on $R$.

We can now argue as for \eqref{stimarI}, where the norms  $\|s\|_{L^2(0,T;L^2(\Gamma))}$ and  $\|c\|_{L^2(0,T;L^2(\Gamma))}$
are controlled by the norms in $\|s\|_{L^2(0,T;H^1(\Omega))}$, $\|c\|_{L^\infty(0,T;H^1(\Omega))}$, respectively.
Thus, we obtain (replacing $\widehat r$ with $r$)
$$\|r\|_{H^1(0,T;L^2(\Gamma))}\leq C,$$
for a constant $C$ independent of $R$.

Finally, recalling \eqref{boundst} and \eqref{stimaPFs4}, we can find
\begin{equation}
\|s\|_{W^{1,\infty}(0,T;L^2(\Omega))\cap H^1(0,T;H^1(\Omega))\cap L^\infty(0,T;H^2(\Omega))}\leq  C,
\end{equation}
for some $C$ independent of $R$.

Choosing now, in the fixed point argument, $R>0$ large enough such that $R \geq T^{1/2} C$, then we can prove that
${\cal S} : {\cal X}_{R,T} \to  {\cal X}_{R,T}$ is a contraction in ${\cal X}_{R,T}$ so that the solution is defined on $(0,T)$.

\section{Numerical Examples}

A fully implicit finite elements scheme has been used to
solve numerically the proposed model where $s, c$ and $r$ are defined within the body or along the boundary $\Gamma$ of the domain $\Omega$.
The evolution equation for $s$ (see \eqref{reazdiff}) is discretized with linear finite elements. The code is based on the freely available Open Source package
deal.II \cite{DealII}.

Computational tests were performed to verify the efficiency
of the numerical technique in 2D setup ($x=x_1 e_1+x_2 e_2$). Moreover, these computations permit to illustrate the main features of the proposed approach.
 The considered domain is a simple marble square domain $1 \times 1$ mm$^2$  invested by a polluted air flow along the left vertical side whereas the other
 faces are isolated. The distribution of the pollutant $SO_2$ on the external border is constant with a concentration equal to $s_0$.
 Moreover, the material is homogeneous so that the initial condition for the calcite density reads $c\left(x,0\right)=c_0$, being $c_0$ a positive constant,
 whereas the concentration of $SO_2$ within the solid is null so $s\left(x,0\right)=0$.

The numerical examples have been obtained by assuming simple expressions for the functions  of \eqref{boundflux}, \eqref{surfacedam}.
In particular, according to the physical evidence that the higher is the rugosity value more exchange surface is available, two monotone relations have been
considered for $\nu(r)$ with $r \geq0$:

\begin{itemize}
  \item linear $\nu \left( r \right)={\nu _0} + \frac{{{\nu _l} - {\nu _0}}}{{{r_l}}}\,r$
  \item parabolic $\nu \left( r \right)={\nu _0} + \frac{{{\nu _l} - {\nu _0}}}{{{r_l}^2}}\,r^2$
\end{itemize}
where the constants are
\begin{itemize}
  \item $\nu _0$: the minimum value for a fully flat surface ($r=0$)
  \item $\nu _l$: the maximum value when it is assumed that $r \in [0,r_l]$ (cf. Remark \ref{stazs3}, setting $r_l = R_0$).
\end{itemize}
The two curves give the same values of $\nu \left( r \right)$ for $r \in \left\{ {0,{r_l}} \right\}$.

For simplicity, no constraints on the variable $r$ are introduced so to have $\partial W(r)=\Psi'(r)=0$ and environmental sources are not considered thus $F=0$.
Under these hypotheses, equation \eqref{surfacedam} reduces to
\begin{equation}\label{reducedr}
\partial_t r+G(r,c,s)= 0,\quad\hbox{ on }\Gamma\times(0,T).
\end{equation}
For the function $G(r,c,s)$ the following expression is adopted
\begin{equation}\label{g_s_c_r}
G(r,c,s)=  - \varphi \left( c \right)cs\left( {{1} + \frac{{r}}{{1 + r}}} \right) g\,,
\end{equation}
where $g$ is a parameter that defines the rate of rugosity evolution and has to be fitted according to experimental evidences.

In the numerical examples different initial rugosity distributions have been considered: piecewise constant or random distribution along the boundary.
In both cases linear and parabolic diffusion functions for $\nu (r)$ have been adopted. Moreover, it has been assumed $\lambda =100$, $g=30$ and the time
is discretized in $n$ constant time intervals $\Delta t=1/5000$.

\subsection{Piecewise constant rugosity}
\begin{figure}[ht]
\centering
\centering
 \begin{tabular}{c c}
 a)\
\includegraphics[height=6.5cm,width=6.5cm]{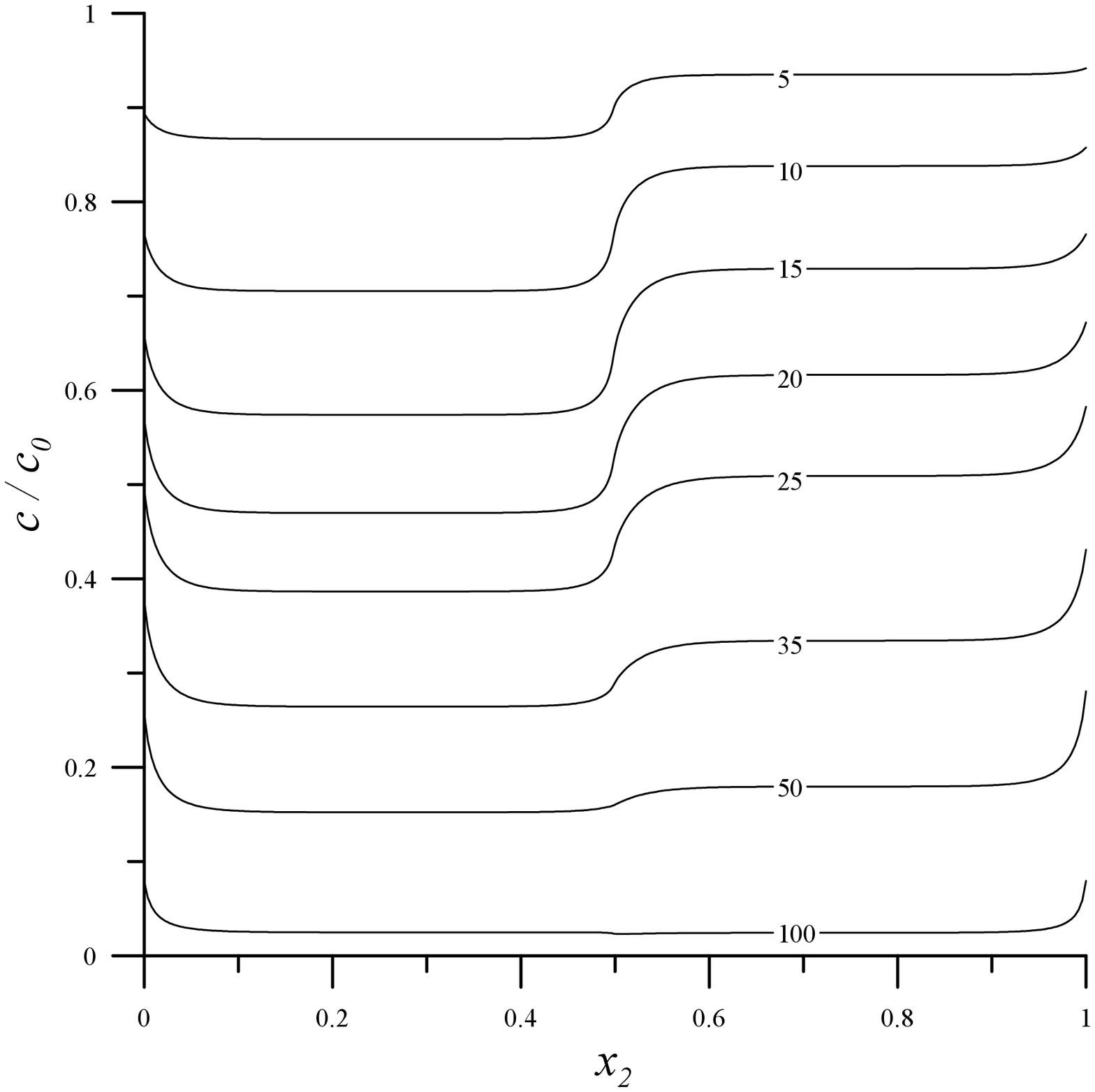} 
  &
  \
  \
  b)\
\includegraphics[height=6.5cm,width=6.5cm]{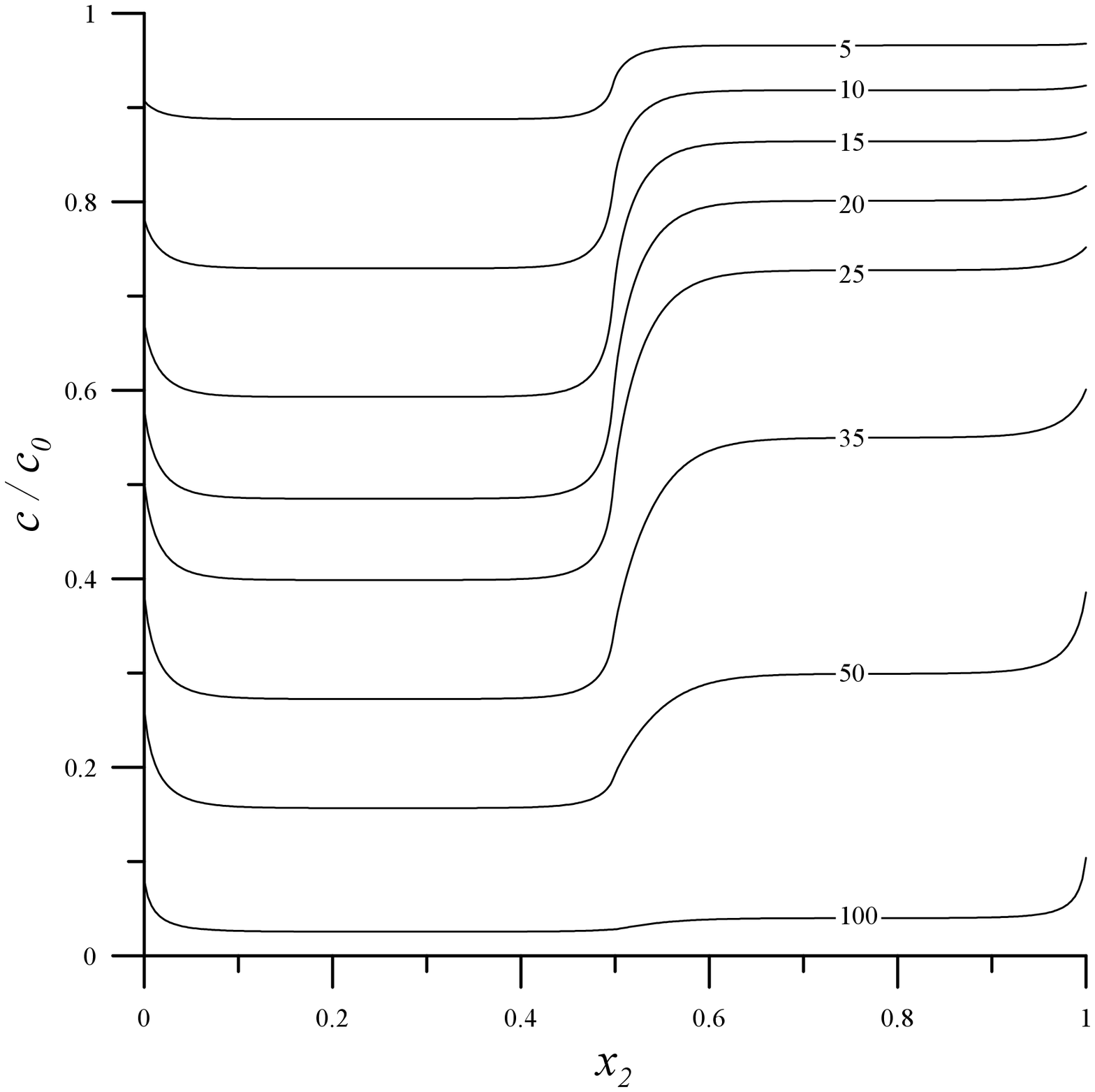} 
  \end{tabular}
 \caption{Evolution of $c$ along the left vertical bounder for $r$ piecewise and $\nu (r)$ a) linear and b) parabolic.} \label{Fig:c_piecewise}
\end{figure}

\begin{figure}[ht]
\centering
\centering
 \begin{tabular}{c c}
 a)\
\includegraphics[height=6.5cm,width=6.5cm]{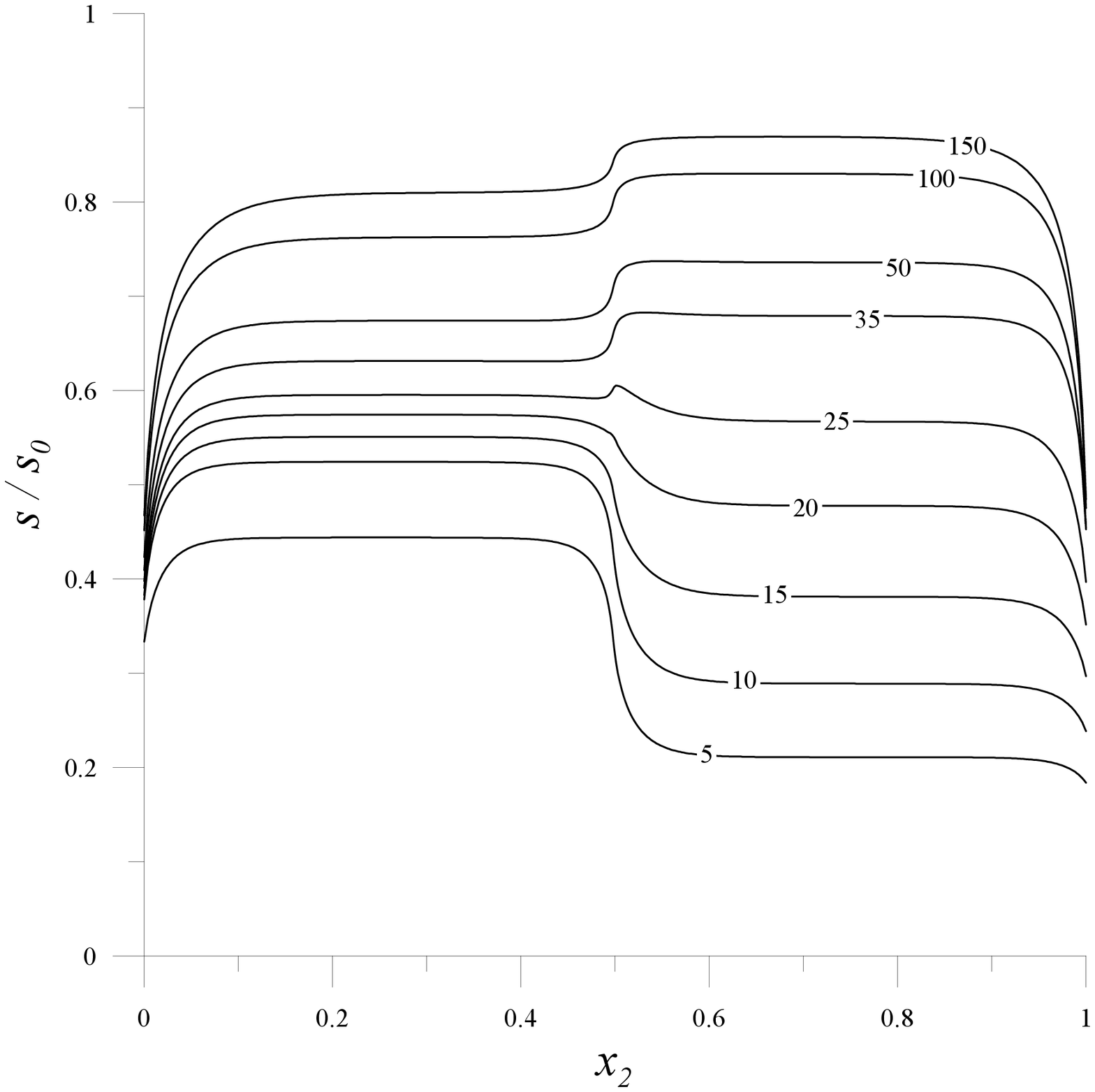} 
  &
  \
  \
  b)\
\includegraphics[height=6.5cm,width=6.5cm]{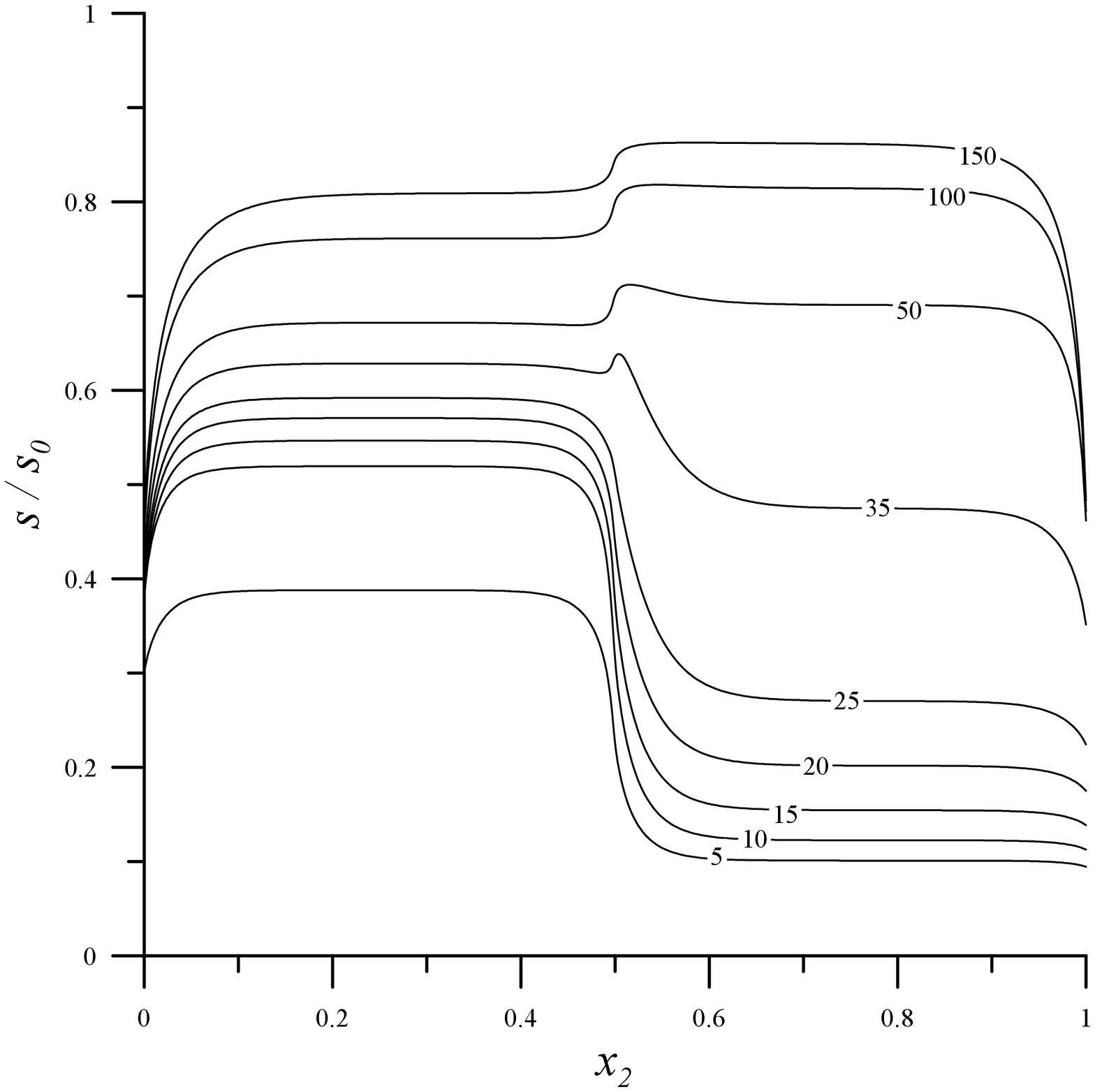} 
  \end{tabular}
 \caption{Evolution of $s$ along the left vertical bounder for $r$ piecewise and $\nu (r)$ a) linear and b) parabolic.} \label{Fig:s_piecewise}
\end{figure}

\begin{figure}[ht]
\centering
\centering
 \begin{tabular}{c c}
 a)\
\includegraphics[height=6.5cm,width=6.5cm]{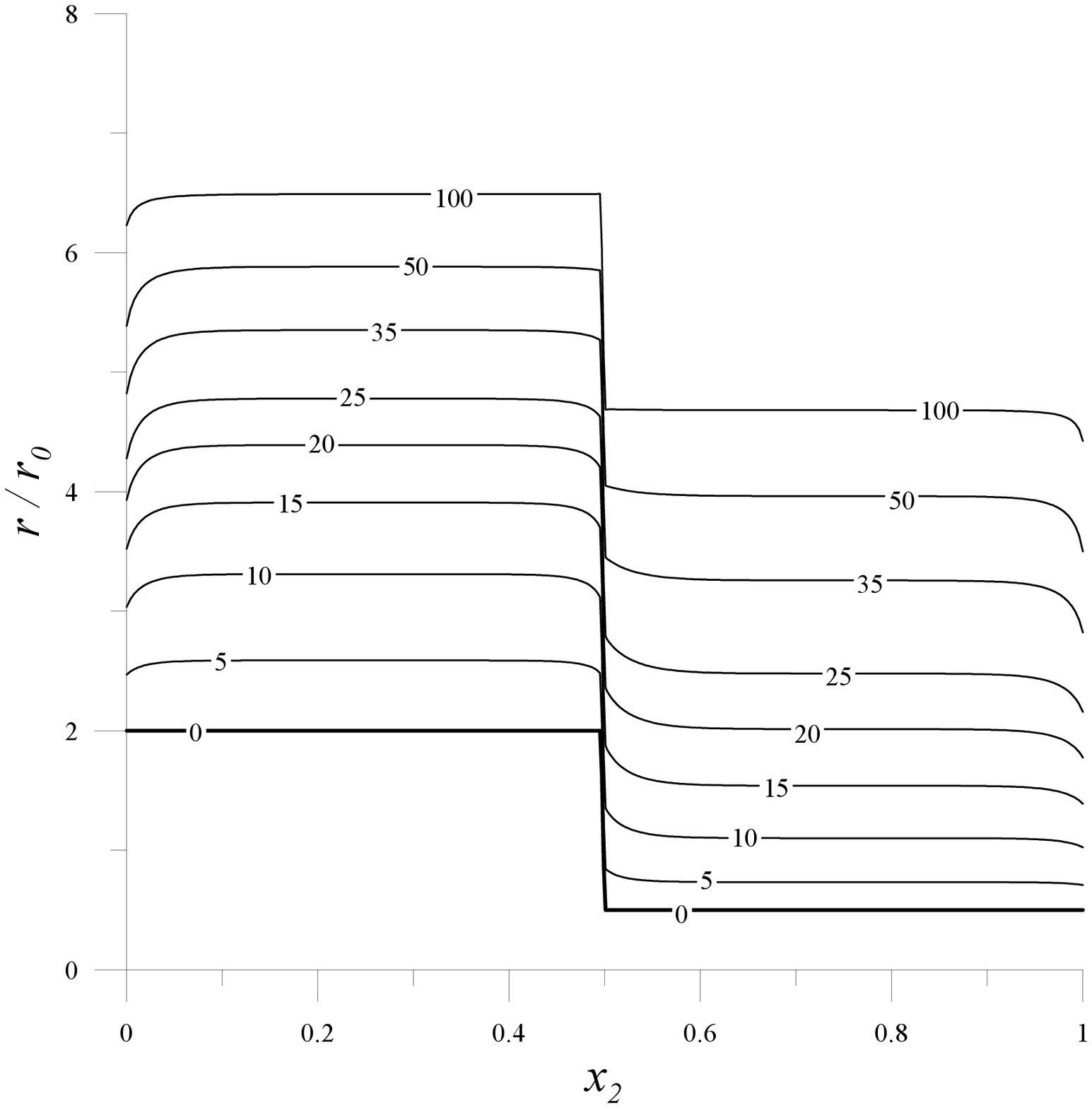} 
  &
  \
  \
  b)\
\includegraphics[height=6.5cm,width=6.5cm]{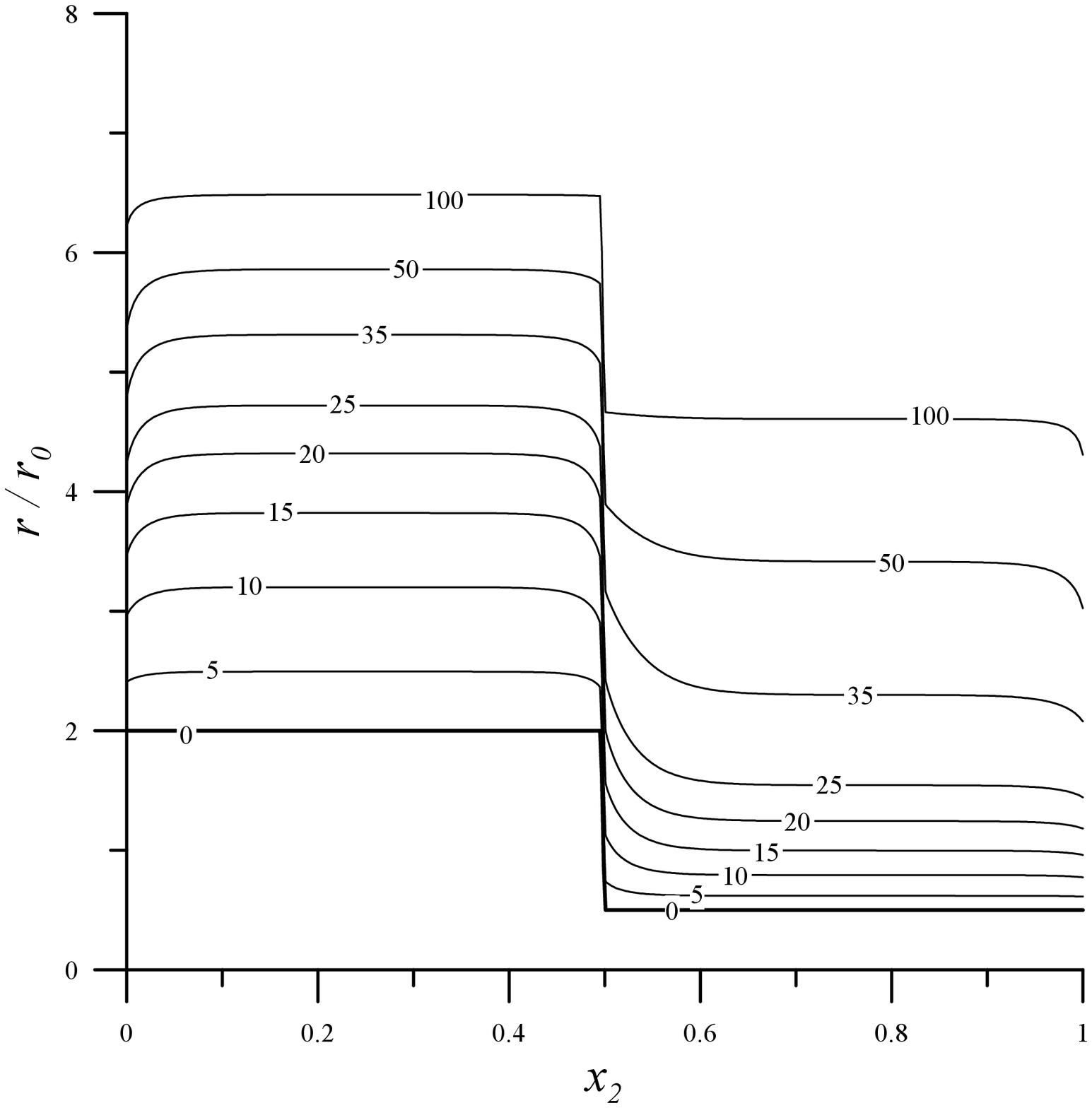} 
  \end{tabular}
 \caption{Evolution of $r$ along the left vertical bounder for $r$ piecewise and $\nu (r)$ a) linear and b) parabolic.} \label{Fig:r_piecewise}
\end{figure}

In the first example, piecewise constant value of the rugosity is assumed along the contour: for $x_2<0.5$, $r\left(x,0\right)=0.5 r_0$
and for $x_2\geq 0.5$, $r\left(x,0\right)=2 r_0$ with $r_0$ a positive constant.
In Figs. \ref{Fig:c_piecewise}, \ref{Fig:s_piecewise}, \ref{Fig:r_piecewise} are reported the profiles of the three variables $c, s, r$
along the vertical side invested by the pollutant $SO_2$ for several time steps. Both linear and parabolic evolution laws for $\nu (r)$ are considered.

The sulphation process reported in Fig. \ref{Fig:c_piecewise}, due to the different assumed values of $r$, evolves with various velocities along the boundary.
In fact, for small time step values, the profile of $c$ presents a significant gradient near $x_2=0.5$ where the initial rugosity is discontinuous.
This phenomenon is much more evident in case the parabolic relation for $\nu (r)$ is assumed.
Subsequently, the solution becomes smoother and tends to be homogenous as $c \rightarrow 0$.
This behavior can be understood analyzing the evolution of $s$ and $r$ as plotted in Figs. \ref{Fig:s_piecewise}, \ref{Fig:r_piecewise} respectively.
In fact, small rugosity values act as a barrier to the penetration of $s$. Even for these two variables high gradients are evident near the middle of the side.
Indeed, the parabolic function for $\nu (r)$ reveals higher gradient values.

\begin{figure}[ht]
\centering
\centering
 \begin{tabular}{c c}
 a)\
\includegraphics[height=6.5cm,width=6.5cm]{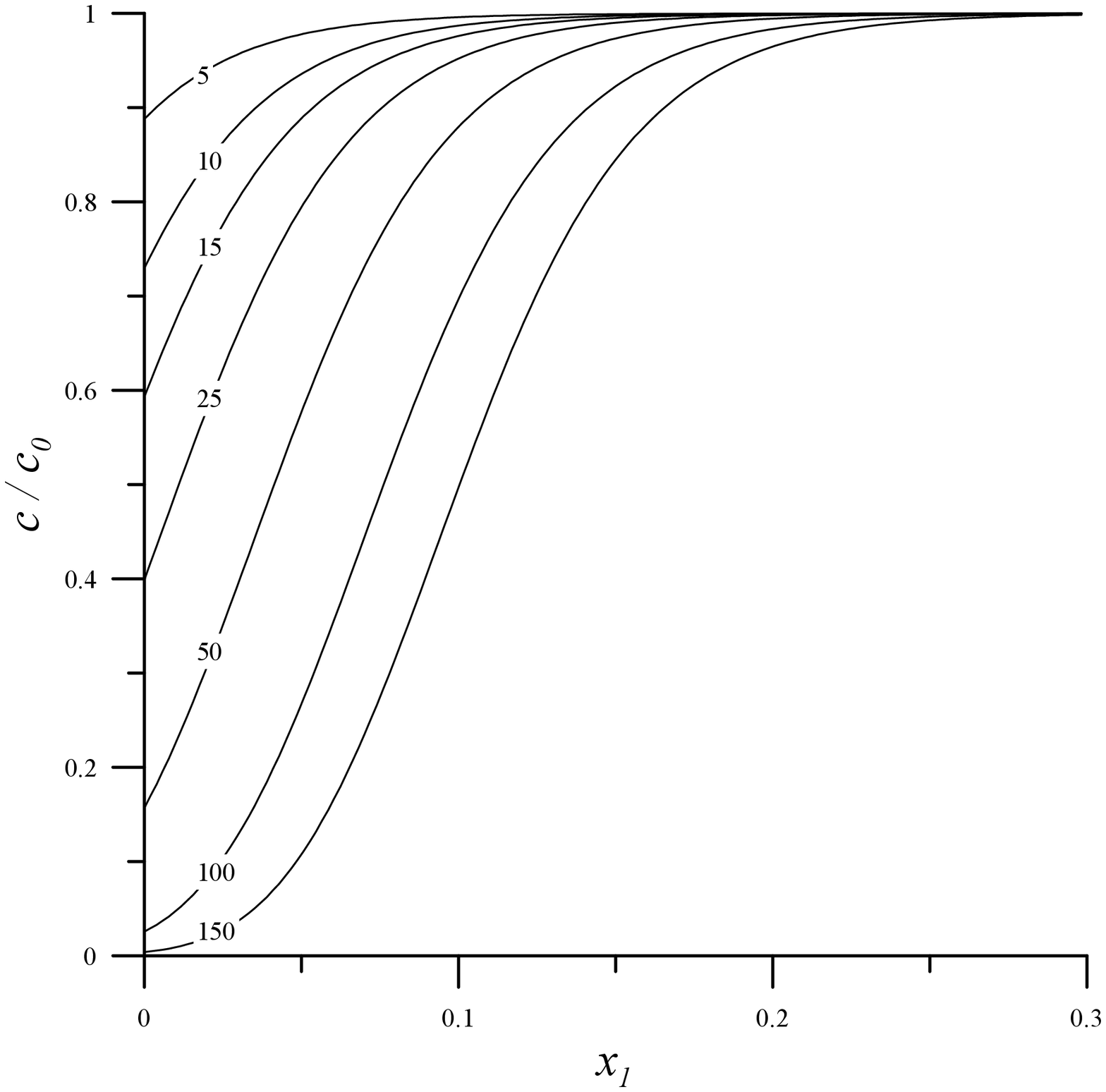} 
  &
  \
  \
  b)\
\includegraphics[height=6.5cm,width=6.5cm]{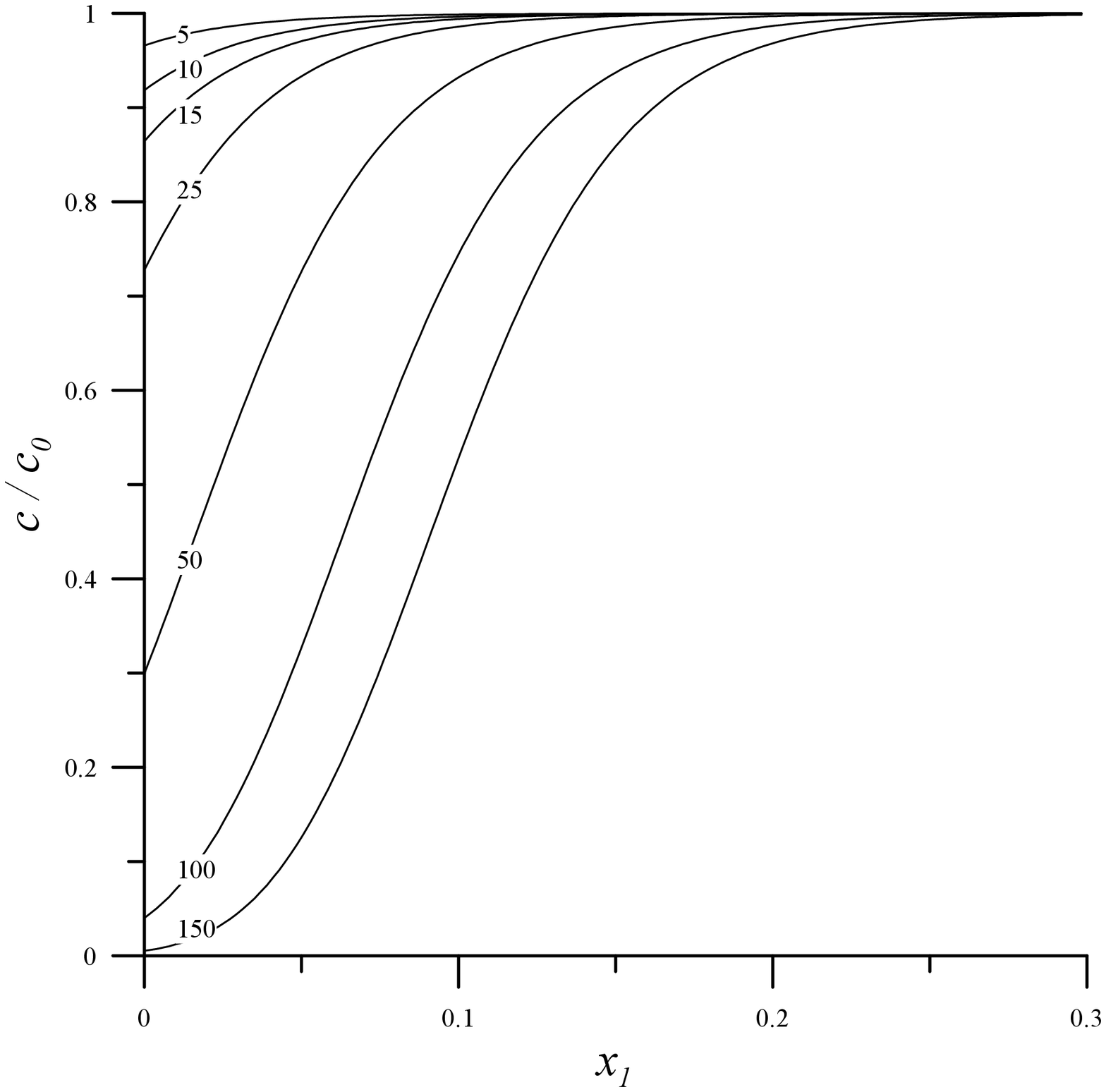} 
  \end{tabular}
 \caption{Evolution of $c$ along a horizontal line within the solid located at a) $x_2=0.25$ and b) $x_2=0.75$
 assuming parabolic relationship for $\nu (r)$.} \label{Fig:c_parabolic_inner}
\end{figure}

\begin{figure}[ht]
\centering
\centering
 \begin{tabular}{c c}
 a)\
\includegraphics[height=6.5cm,width=6.5cm]{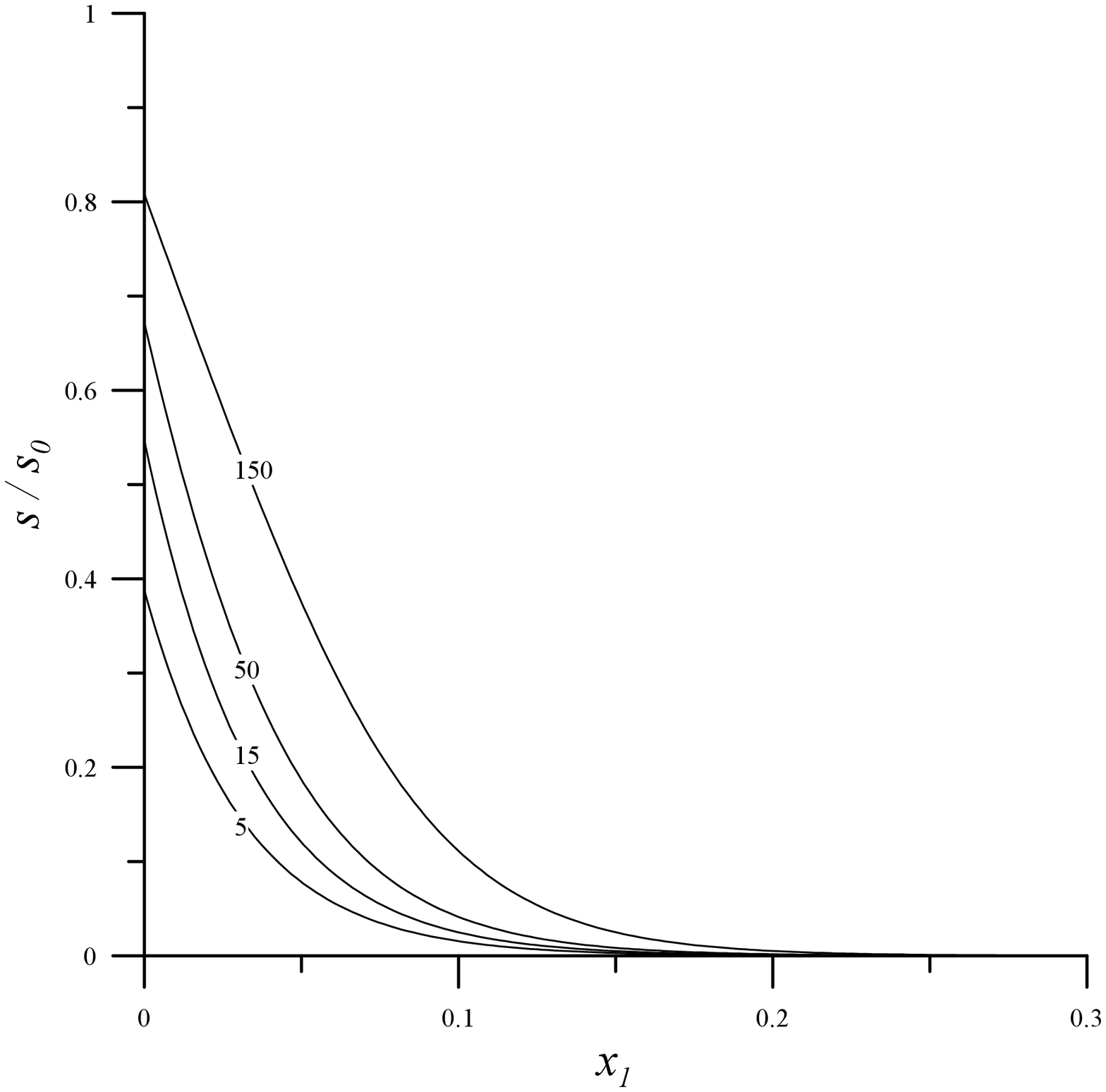} 
  &
  \
  \
  b)\
\includegraphics[height=6.5cm,width=6.5cm]{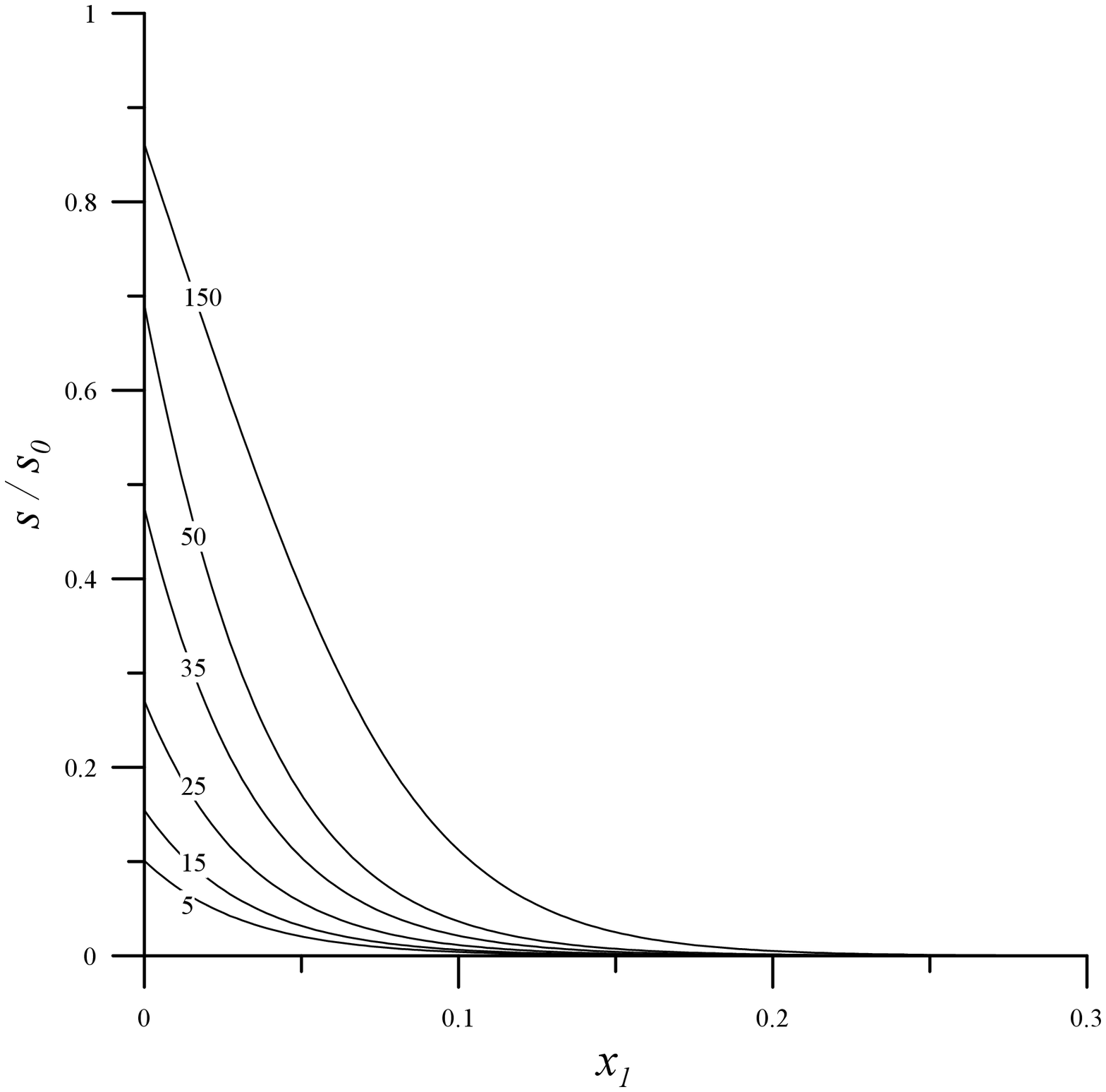} 
  \end{tabular}
 \caption{Evolution of $s$  along a horizontal line within the solid located at a) $x_2=0.25$ and b) $x_2=0.75$
 assuming parabolic relationship for $\nu (r)$.} \label{Fig:s_parabolic_inner}
\end{figure}

\begin{figure}[ht]
\centering
\centering
 \begin{tabular}{c c}
 a)\
\includegraphics[height=5.5cm,width=6.5cm]{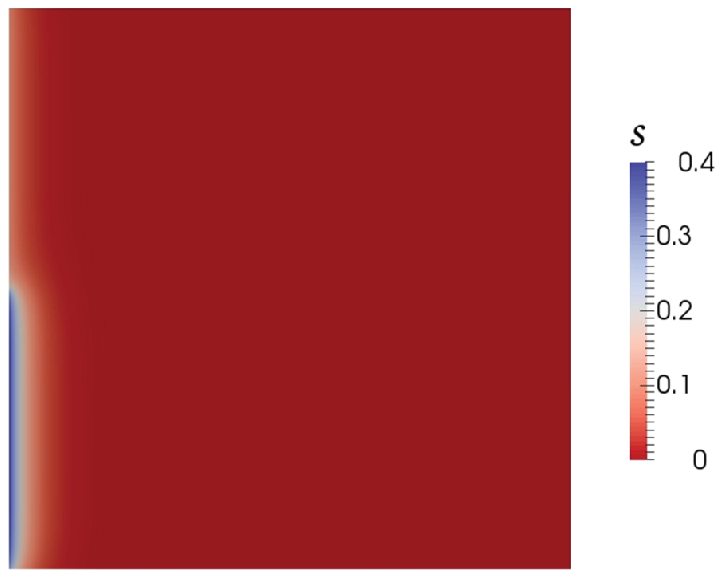} 
  &
  \
  \
  b)\
\includegraphics[height=5.5cm,width=6.5cm]{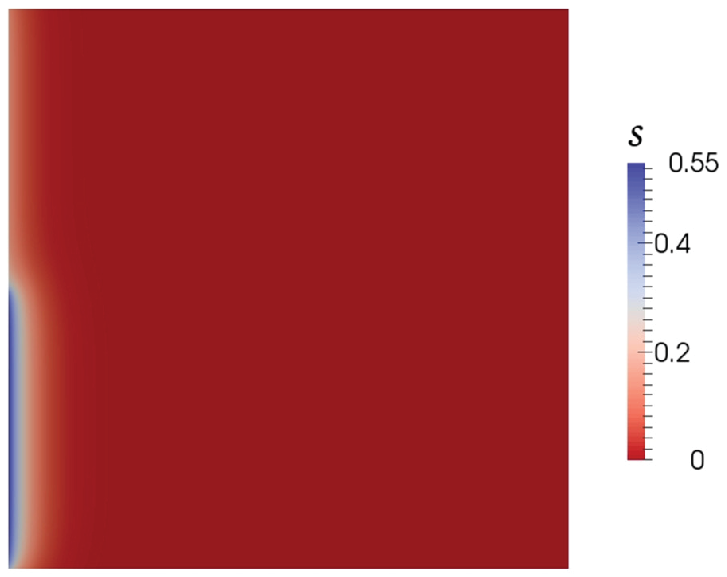} 
  \end{tabular}
 \caption{Concentration of $S0_2$  within the solid at different time step a) $n=5$ and b) $n=15$
 assuming parabolic relationship for $\nu (r)$.} \label{Fig:s_maps_parabolic_inner}
\end{figure}

The effect of different rugosity values also influences the solution within the solids. The evolutions of the variables $c$ and $s$ are plotted along two
horizontal lines located at $x_2=\{0.25, 0.75\}$ for different time steps in Figs. \ref{Fig:c_parabolic_inner}, \ref{Fig:s_parabolic_inner}.
At $x_2=0.25$, where the initial rugosity is high, the sulphation process is affected by the rapid diffusion of $s$ within the solid that activates the
calcite transformation. On the contrary, the diffusion of $SO_2$ at $x_2=0.75$, and consequently the marble degradation, begins slowly due to the initial
small value of $\nu (r)$. For large time ($n>100$) the influence of the surface rugosity turns to be negligible on the inner solutions.
The maps of $s$ within the solid is plotted in Fig. \ref{Fig:s_maps_parabolic_inner} for time steps $n=\{5, 15\}$ and parabolic $\nu (r)$.
In the lower portion of the domain, due to higher rugosity values along the border, larger diffusion of $SO_2$ is evident.

\FloatBarrier

\subsection{Random rugosity}
\begin{figure}[ht]
\centering
\centering
 \begin{tabular}{c c}
 a)\
\includegraphics[height=6.5cm,width=6.5cm]{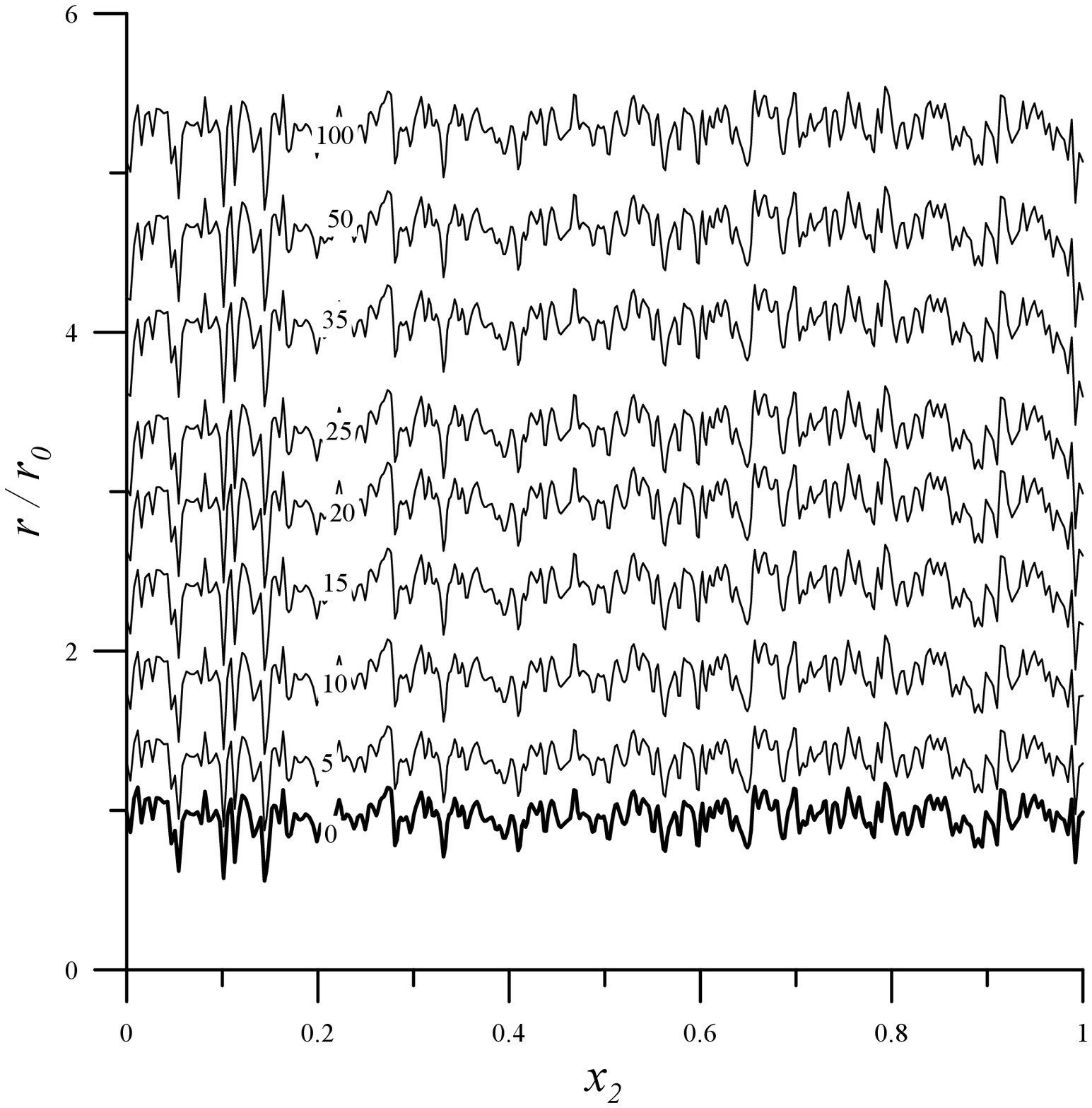} 
  &
  \
  \
  b)\
\includegraphics[height=6.5cm,width=6.5cm]{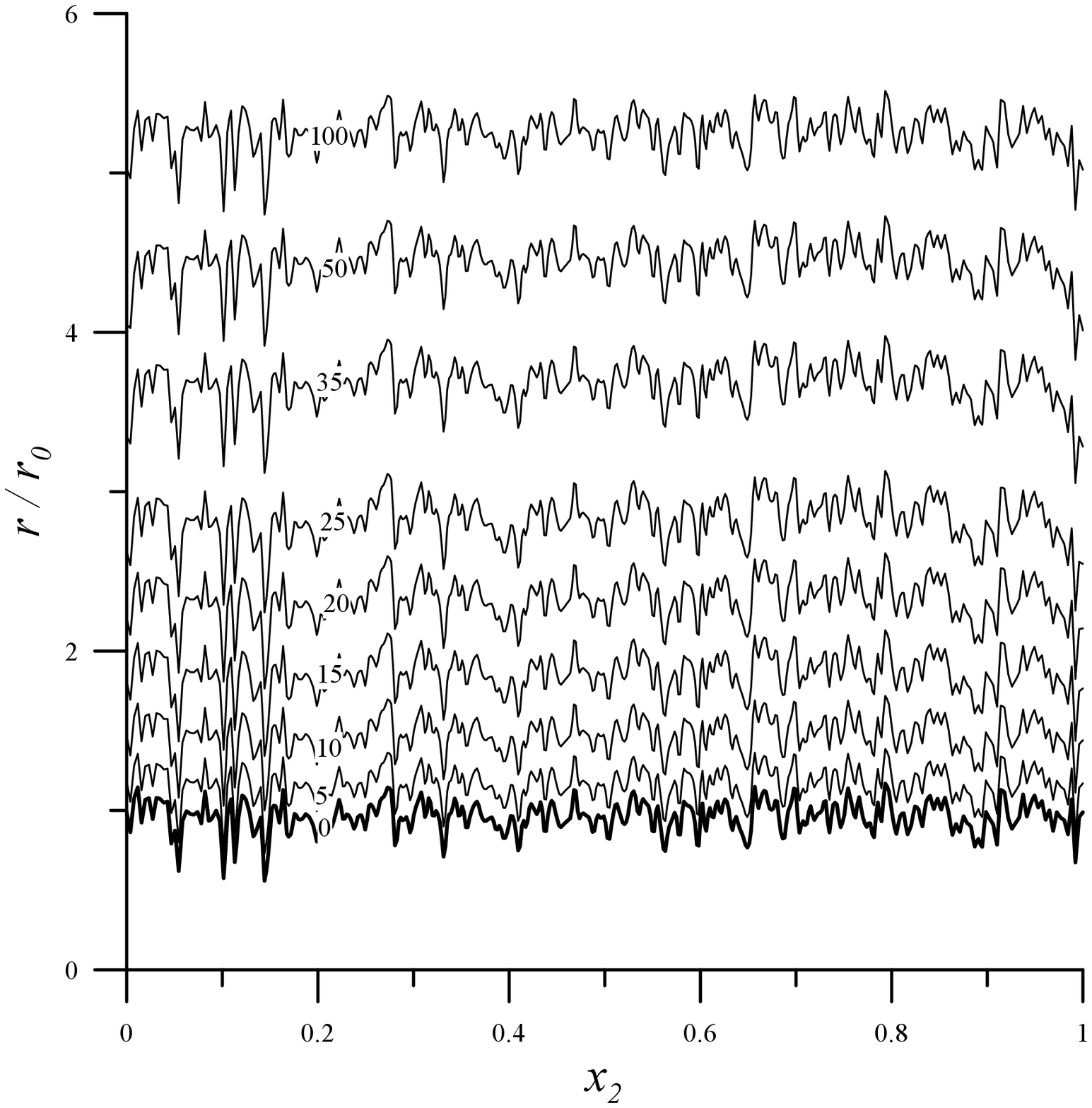} 
  \end{tabular}
 \caption{Evolution of $r$ along the left vertical bounder for $r$ random and $\nu (r)$ a) linear and b) parabolic.} \label{Fig:r_random}
\end{figure}

Subsequently, a random surface rugosity based on Weibull's statistics is assigned on $\Gamma$.
The initial value of $\bar r_0$ is therefore assumed as
\begin{equation}
	\bar r_0 =  r_0\left( { \ln{\dfrac{1}{1-\lambda}}}\right) ^{1/m}
\label{random_r_0}
\end{equation}
with $r_0$ and $m$ being, respectively, the Weibull shape and modulus parameters, and $\lambda$ a random variable ranging from~0 to~1.
Here, $ r_0$ is the mean value of surface rugosity. In the computation $m=10$ has been considered.
The initial rugosity profile is reported in Fig. \ref{Fig:r_random}. In the evolution, the asperities are maintained for $\nu (r)$
linear whereas are emphasized in case of parabolic relation. For large time the rugosity does not evolve anymore and the heterogenous profile is preserved.

\begin{figure}[ht]
\centering
\centering
 \begin{tabular}{c c}
 a)\
\includegraphics[height=6.5cm,width=6.5cm]{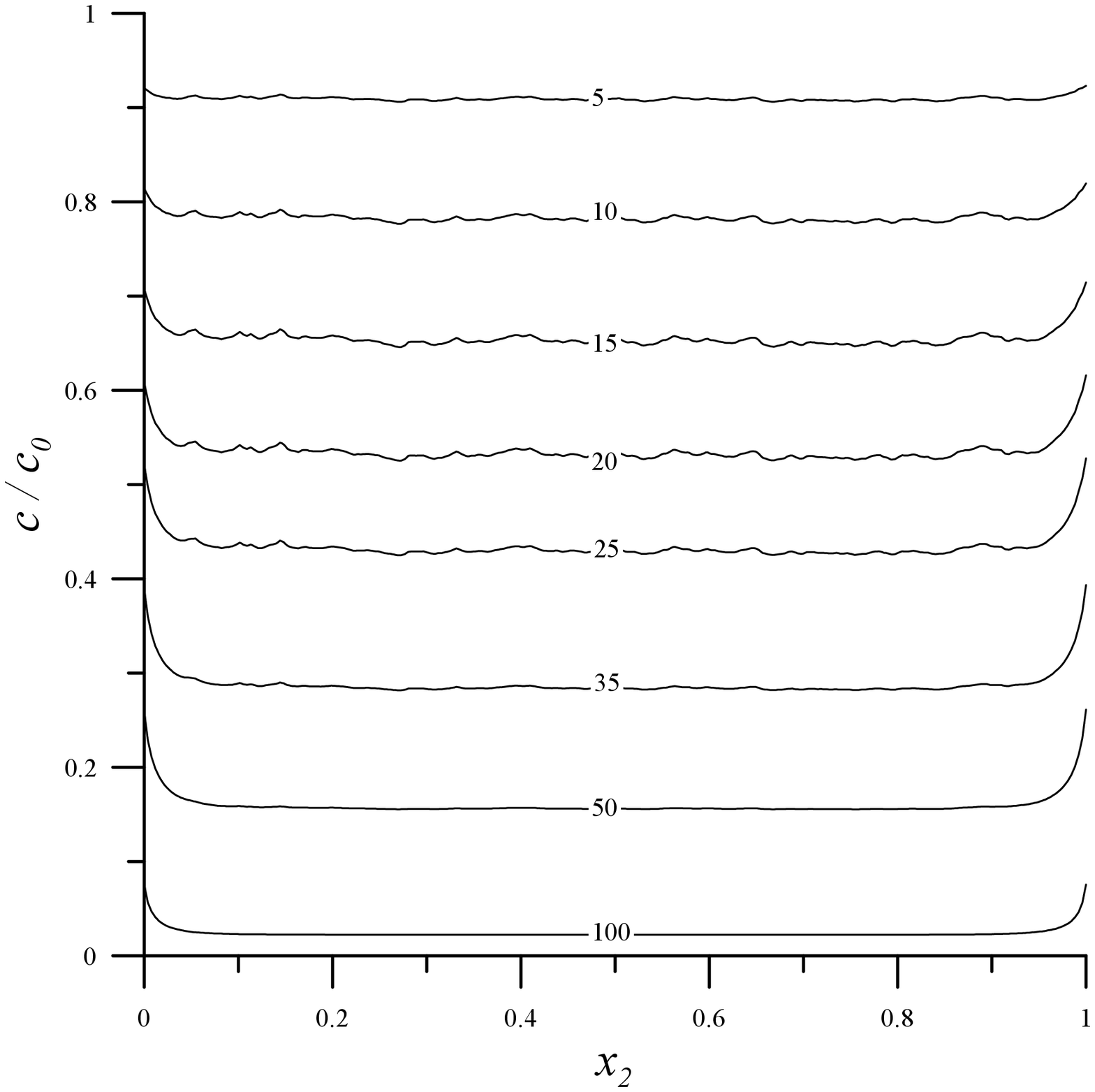} 
  &
  \
  \
  b)\
\includegraphics[height=6.5cm,width=6.5cm]{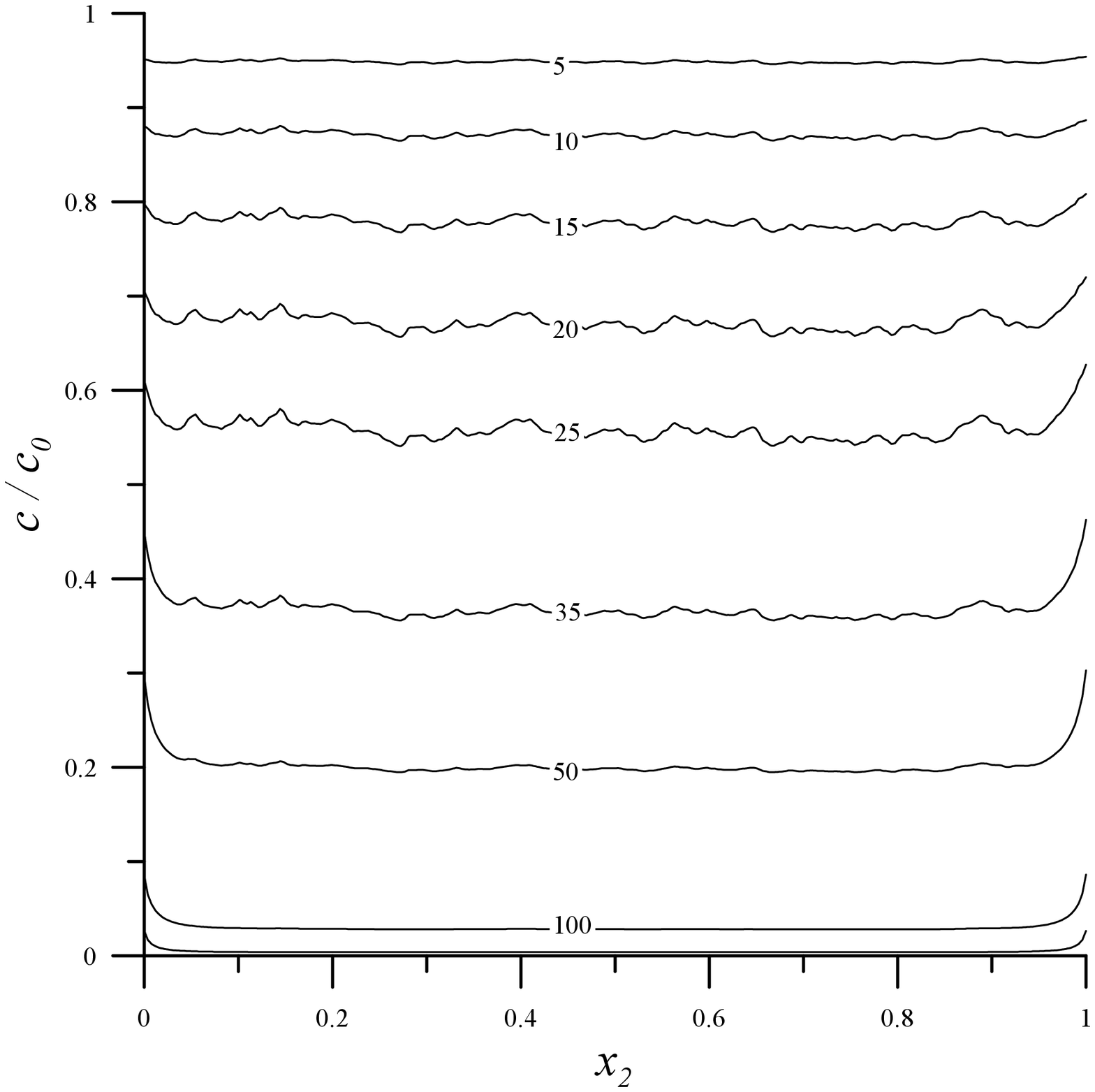} 
  \end{tabular}
 \caption{Evolution of $c$ along the left vertical bounder for $r$ random and $\nu (r)$ a) linear and b) parabolic.} \label{Fig:c_random}
\end{figure}

\begin{figure}[ht]
\centering
\centering
 \begin{tabular}{c c}
 a)\
\includegraphics[height=6.5cm,width=6.5cm]{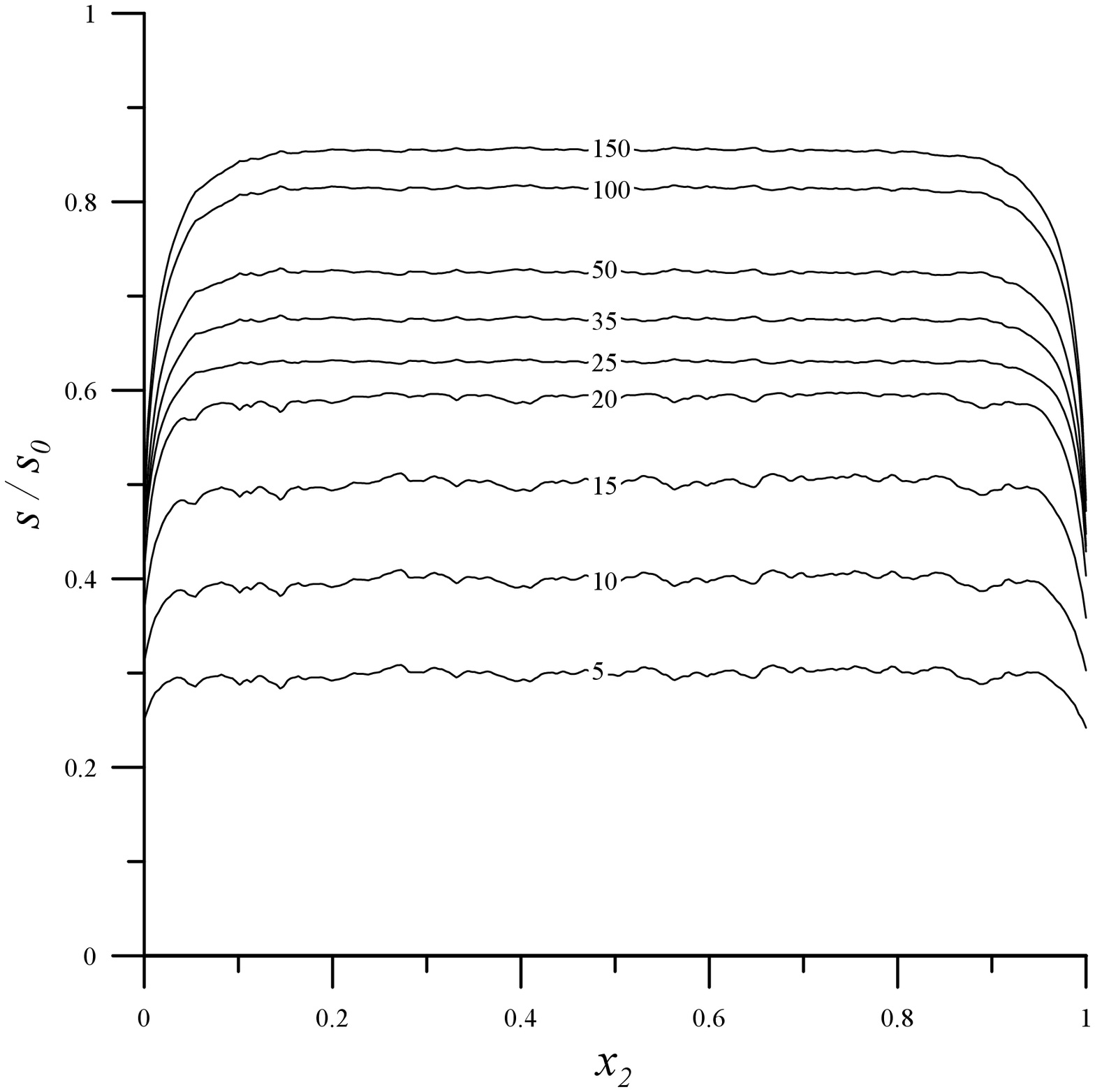} 
  &
  \
  \
  b)\
\includegraphics[height=6.5cm,width=6.5cm]{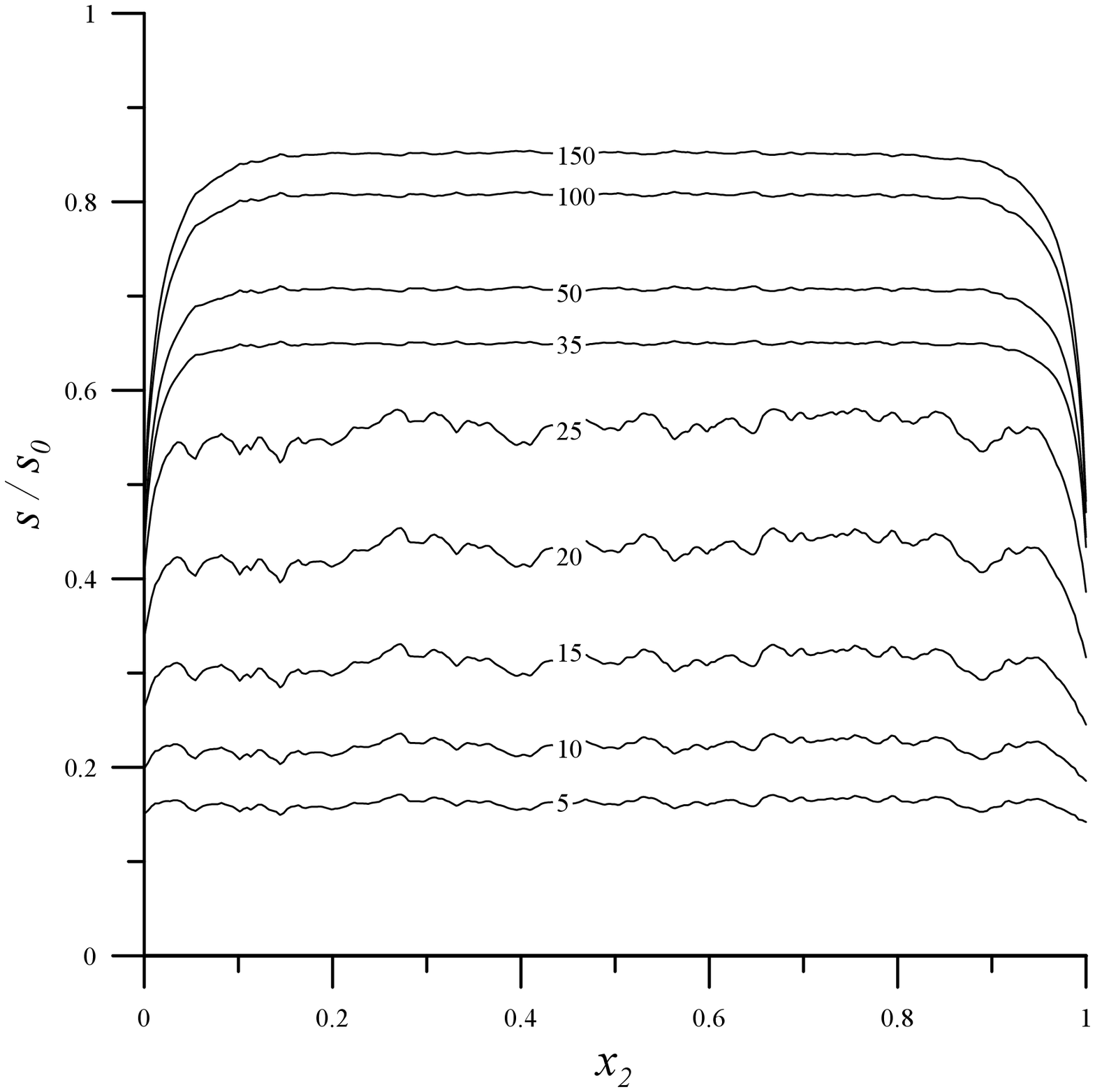} 
  \end{tabular}
 \caption{Evolution of $s$ along the left vertical bounder for $r$ random and $\nu (r)$ a) linear and b) parabolic.} \label{Fig:s_random}
\end{figure}

The initial homogenous distribution of $c$ is abandoned during the sulphation process as clearly stated by Fig. \ref{Fig:c_random} that
reports the evolution of $c$ for $\nu (r)$ linear and parabolic. The profile becomes extremely jagged for parabolic relationship of $\nu (r)$.
Subsequently, the asperities vanish as $c \rightarrow 0$. Analogous behavior appears for the $SO_2$ concentration as reported in Fig. \ref{Fig:s_random}.
The rugosity profile, as state in Fig. \ref{Fig:r_random}, maintains the shape with an increment of the fluctuation.

\FloatBarrier

\section{Conclusions}

{\elefin The local in time well-posedness and  the global existence theorem are preliminary but essential results related to our problem.
A nontrivial issue might be the analysis of the longtime behavior of solutions and, in particular, their possible convergence to a steady state. Moreover,
it would be interesting (and challenging) to extend our theoretical analysis to significant generalizations which account for further features of this
complex phenomenon.  In particular, in the coupling of bulk-surface equations, further actions should be taken into account, like the bulk damage and
thermo-mechanical effects forcing the surface deterioration.}

The presented numerical simulations highlight the main features of the proposed model.
The introduction of the superficial rugosity permits to account for an important physical property of the material that strongly influences
initiation and the evolution of the sulphation stone subjected to pollutant. In fact, high rugosity values offers wider surface for the $SO_2$ diffusion that
induces the sulphation process. On the other hand, extremely smooth surfaces slow down the degradation process.

At the present stage, two natural development are possible. Firstly, an extensive experimental campaign has to be performed in order to fit the parameters of
the proposed model and to define a proper function for $\nu (r)$. Secondly, the fracture phenomenon and mechanical degradation has to be coupled with the
process of sulphation of the stones. In particular, the best candidate is the variational formulation of fracture mechanics that recently has gained much
attention and demonstrated to be effective to  determine complex crack paths \cite{F1}.

\section{Ackowledgements}

The work of C. Cavaterra was supported by the FP7-IDEAS-ERC-StG 256872 (EntroPhase). E. Bonetti, C. Cavaterra and M. Grasselli are members of
the Gruppo Nazionale per l'Analisi Matematica, la Probabilit\`{a} e le loro Applicazioni (GNAMPA) of the Istituto Nazionale di Alta Matematica (INdAM).


\begin{thebibliography}{99}

\bibitem{N1}
Aregba-Driollet D.,  Diele  F., Natalini, R.,
A mathematical model for the sulphur dioxide aggression to calcium carbonate stones: numerical approximation and asymptotic analysis,
{\sl SIAM J. Appl. Math.},  {\bf 64} (2004), 1636-1667.



\bibitem{baiocchi}
Baiocchi C., Sulle equazioni differenziali astratte lineari del primo e del secondo ordine negli spazi di Hilbert,
{\sl Ann. Mat. Pura Appl.}, (4) {\bf 76} (1967), 233-304.


\bibitem{DealII}
Bangerth W., Hartmann R., Kanschat G., {deal.II} -- a General Purpose Object Oriented Finite Element Library, {\sl ACM Trans. Math. Softw.}, {\bf 33} (2007), 1-27.


\bibitem{Bonetti}
Bonetti E., Fr\'emond M.,  Analytical results on a model for damaging in domains and interfaces, {\sl ESAIM Control Optim. Calc. Var.}, {\bf 17} (2011), 955-974.

\bibitem{N2}
Clarelli F., Fasano A., Natalini R., Mathematics and monument conservation: free boundary models of marble sulfation,
{\sl SIAM J. Appl. Math.}, {\bf  69} (2008), 149-168.

\bibitem{F1}
Freddi F., Royer-Carfagni G., Regularized Variational Theories of Fracture: a Unified Approach, {\sl Journal of the Mechanics and Physics of Solids},
{\bf 58} (2010), 1154-1174.



\bibitem{Fremond}
Fr\'emond M., Non-smooth thermomechanics. Springer-Verlag, Berlin, 2002.


\bibitem{N3}
Giavarini C., Santarelli M.L., Natalini R., Freddi F., A non-linear model of sulphation of porous stones: Numerical simulations and preliminary
laboratory assessments, {\sl Journal of Cultural Heritage}, {\bf 9} (2008), 14-22.


\bibitem{GN1}
Guarguaglini F.R, Natalini R., {Global existence of solutions to a nonlinear model of sulphation phenomena in calcium carbonate stones},
{\sl Nonlinear Anal}, {\bf 6} (2005), 477-494.

\bibitem{GN3D}
Guarguaglini F.R, Natalini R., Nonlinear transmission problems for quasilinear diffusion systems,  {\sl Netw. Heterog. Media}, {\bf 2} (2007), 359-381.

\bibitem{GN2}
Guarguaglini F.R., Natalini R., Fast reaction limit and large time behavior of solutions to a nonlinear model for sulphation phenomena,
{\sl Comm. Partial Differential Equations}, {\bf 32} (2007), 163-189.

\bibitem{peletier}
 Mielke A., A gradient structure for reaction-diffusion systems and for energy-drift-diffusion systems, {\sl Nonlinearity}, {\bf  24} (2011), 1329-1346.


\bibitem{Natalini-BARENBLAT}
 Natalini R.,  Nitsch C., Pontrelli G., Sbaraglia S., A numerical study of a nonlocal model of damage propagation under chemical aggression,
 {\sl European J. Appl. Math.}, {\bf 14 } (2003), 447-464.


\bibitem{pao}
Pao C. V., Nonlinear parabolic and elliptic equations. Plenum Press, New York, 1992.



\bibitem{Taylor}
Taylor M., Partial Differential equations I. Basic theory. Applied Mathematical Sciences, {\bf 115},  Springer, New York, 2011.


\bibitem{whitehouse}
Whitehouse D., Surfaces and their Measurement. Butterworth-Heinemann, Boston, 2012.

\bibitem{iso}
BS EN ISO 4287:2000, Geometrical product specification (GPS). Surface texture. Profile method. Terms, definitions and surface texture parameters.



\end{thebibliography}
\end{document}